\documentclass{amsart}
\usepackage{amsfonts,amscd,amsmath,amssymb,amsthm,mathrsfs}
\usepackage{array,latexsym,color,a4wide,float,enumerate}
\usepackage{graphicx,epsfig}
\usepackage[all]{xy}
\usepackage{color,mathrsfs}
\usepackage{longtable}
\usepackage[hypertex,linkbordercolor={0 0 1}]{hyperref}

\theoremstyle{plain}
\newtheorem{theorem}{Theorem}[section]
\newtheorem{lemma}[theorem]{Lemma}
\newtheorem{proposition}[theorem]{Proposition}
\newtheorem{corrolary}[theorem]{Corollary}

\theoremstyle{definition}
\newtheorem{definition}[theorem]{Definition}
\newtheorem{example}[theorem]{Example}
\newtheorem{noname}[theorem]{}
\newtheorem{remark}[theorem]{Remark}
\newtheorem{construction}[theorem]{Construction}
\newtheorem{notation}[theorem]{Notation}

\theoremstyle{remark}
\newtheorem*{smallremark}{Remark}

\newtheorem{case}{Case} \makeatletter \@addtoreset{case}{theorem}\makeatother
\newtheorem{claim}{Claim} \makeatletter \@addtoreset{claim}{theorem}\makeatother

\newcommand{\bthm}{\begin{theorem}}
\newcommand{\bprop}{\begin{proposition}}
\newcommand{\blem}{\begin{lemma}}
\newcommand{\bcor}{\begin{corrolary}}
\newcommand{\brem}{\begin{remark}}
\newcommand{\bdfn}{\begin{definition}}
\newcommand{\bitem}{\begin{itemize}}
\newcommand{\bex}{\begin{example}}
\newcommand{\bno}{\begin{noname}}
\newcommand{\bsrem}{\begin{smallremark}}
\newcommand{\bnot}{\begin{notation}}
\newcommand{\bcon}{\begin{construction}}
\newcommand{\bca}{\begin{case}}
\newcommand{\bcl}{\begin{claim}}

\newcommand{\ecl}{\end{claim}}
\newcommand{\eca}{\end{case}}
\newcommand{\econ}{\end{construction}}
\newcommand{\enot}{\end{notation}}
\newcommand{\esrem}{\end{smallremark}}
\newcommand{\eno}{\end{noname}}
\newcommand{\eex}{\end{example}}
\newcommand{\eitem}{\end{itemize}}
\newcommand{\ethm}{\end{theorem}}
\newcommand{\eprop}{\end{proposition}}
\newcommand{\elem}{\end{lemma}}
\newcommand{\ecor}{\end{corrolary}}
\newcommand{\erem}{\end{remark}}
\newcommand{\edfn}{\end{definition}}

\newcommand{\ba}{\begin{array}}
\newcommand{\ea}{\end{array}}
\newcommand{\benum}{\begin{enumerate}}
\newcommand{\eenum}{\end{enumerate}}

\newcommand{\un}{\underline}
\newcommand{\ov}{\overline}

\newcommand{\wt}{\widetilde}

\def\iso{\cong}

\def\8{\infty}
\def\.{\cdot}
\def\PP{\mathbb{P}}

\def\C{\mathbb{C}}
\def\Z{\mathbb{Z}}
\def\N{\mathbb{N}}
\def\Q{\mathbb{Q}}
\def\R{\mathbb{R}}
\def\E{\widehat{E}}
\def\ovk{\overline\kappa}
\def\kp{\kappa}

\def\:{\colon}

\def\med{\medskip}
\def\ssk{\smallskip}
\def\bsk{\bigskip}
\newcommand{\noin}{\noindent}
\def\Bk{\operatorname{Bk}}

\def\Aut{\operatorname{Aut}}

\def\Sing{\operatorname{Sing}}

\def\Supp{\operatorname{Supp}}

\def\dim{\operatorname{dim}}

\begin{document}

\title[Exceptional singular $\Q$-homology planes]{Exceptional singular $\Q$-homology planes}
\author[Karol Palka]{Karol Palka}
\address{Karol Palka: Institute of Mathematics, University of Warsaw, ul. Banacha 2, 02-097 Warsaw, Poland}
\address{Institute of Mathematics, Polish Academy of Sciences, ul. \'{S}niadeckich 8, 00-956 Warsaw, Poland}
\thanks{The author was supported by Polish Grant MNiSzW (N N201 2653 33)}
\email{palka@impan.pl}
\subjclass[2000]{Primary: 14R05; Secondary: 14J17, 14J26}
\keywords{Acyclic surface, homology plane, exceptional Q-homology plane}

\begin{abstract} We consider singular $\Q$-acyclic surfaces with smooth locus of non-general type. We prove that if the singularities are topologically rational then the smooth locus is $\C^1$- or $\C^*$-ruled or the surface is up to isomorphism one of two exceptional surfaces of Kodaira dimension zero. For both exceptional surfaces the Kodaira dimension of the smooth locus is zero and the singular locus consists of a unique point of type $A_1$ and $A_2$ respectively. \end{abstract}

\maketitle

We consider complex algebraic varieties.

\section{Main result}\label{s:main result}

Because of their homological similarity to $\C^2$ smooth $\Q$-acyclic surfaces serve as a class of test examples for working hypotheses as well for conjectures like cancellation problem or the Jacobian Conjecture, they appear naturally also when studying exotic structures on $\C^n$'s (see \cite[\S 3.4]{Miyan-OpenSurf} for what is known about them).

\bdfn A surface is a \emph{$\Q$-homology plane} if it is normal and $\Q$-acyclic, i.e. $H^*(-,\Q)\cong \Q$.\edfn

A singular $\Q$-homology plane is \emph{logarithmic} if and only if it has at most quotient singularities, i.e. analytically it is locally of type $\C^2/G$ for some finite subgroup $G<GL(2,\C)$. Note that logarithmic $\Q$-homology planes are rational by \cite{GP,PS-rationality}. Singular $\Q$-homology planes appear for example as quotients of smooth ones by the actions of finite groups or as two-dimensional quotients of $\C^n$ by the actions of reductive groups (cf. \cite{KR-ContrSurf}, \cite{Gurjar-quotients}). Let $S'$ be a $\Q$-homology plane and let $S_0$ be its smooth locus ($S'=S_0$ if $S'$ is smooth). Assume that $S_0$ is not of general type, i.e. its Kodaira dimension $\ovk(S_0)$ is smaller than two. The description of these surfaces divides into three main cases depending on the properties of $S_0$: (a) $S_0$ is $\C^1$-ruled, (b) $S_0$ is $\C^*$-ruled, (c) $S_0$ is neither $\C^1$- nor $\C^*$-ruled.

\bdfn A $\Q$-homology plane whose smooth locus is not of general type and is neither $\C^1$- nor $\C^*$-ruled is \emph{exceptional}.\edfn

For non-exceptional $\Q$-homology planes the analysis reduces to the description of singular fibers of respective rulings using the $\Q$-acyclicity. Case (a) and part of case (b) (when $S'$ is logarithmic and the $\C^*$-ruling of $S_0$ extends to a $\C^*$-ruling of $S'$) have been done in \cite{MiSu-hPlanes}. The precise classification and the rest of part (b) will be done in our forthcoming paper. By general structure theorems for open surfaces an exceptional $\Q$-homology plane necessarily has $\ovk(S_0)=0$ (cf. \cite[2.1.1]{Miyan-OpenSurf}, \cite[2.3]{Kawamata2}). The description of smooth exceptional $\Q$-homology planes can be found in \cite[\S 8]{Fujita}. The classification of non-smooth exceptional $\Q$-homology planes is the main goal of this paper. We will do this under some mild assumption on singularities. 

\bdfn A singular point on a normal surface is a \emph{topologically rational singularity} if and only if there exists a resolution of this surface with a rational tree as an exceptional divisor.\edfn

Notice that the singularity is topologically rational if and only if it is \emph{quasirational} (cf. \cite{Ab}) and the dual graph of the respective exceptional locus contains no loops. The class of topologically rational singularities includes the class of rational singularities and is much broader than the class of the quotient ones. Our main result is:

\bthm\label{thm:main result} Up to isomorphism there are exactly two exceptional singular $\Q$-homology planes with at most topologically rational singularities. Both have Kodaira dimension zero and have unique singular points of type $A_1$ and $A_2$ respectively.\ethm

One of the above surfaces comes from the famous \emph{dual Hesse configuration} $(12_3,9_4)$ of points and lines on $\PP^2$ not realizable in $\R\PP^2$ and the second one from the \emph{complete quadrangle} $(4_3, 6_2)$ (see \ref{rem:configurations} and \ref{lem:Y2c_contraction to lines}). We want to emphasize that having some topological results about general singular $\Q$-homology planes (which we obtain in a forthcoming paper) one can easily show that in the above situation the assumption about topological rationality can be omitted with no change for the thesis. However, it is not true that all singular $\Q$-homology planes have topologically rational singularities.

As for now there is no description of $\Q$-homology planes with smooth locus of general type. There are some partial results (see \cite{DieckPetrie}, \cite{Zaidenberg-no_affine_lines, Zaidenberg-no_affine_lines_correction}, \cite{MiTs-absence_of_affine_lines}, \cite{GM-Affine-lines}, \cite{KR-ContrSurf}).

The outline of the proof of the theorem is as follows. First with the help of Bogomolov-Miyaoka-Yau inequality we show in section \ref{s:basic_properties} that each smooth rational curve contained in the snc-minimal smooth completion of $S_0$ has at least two common points with some connected component of the boundary (i.e. it is not \emph{simple}), which in particular shows that $S_0$ is minimal in the open sense (see \cite[2.3.11]{Miyan-OpenSurf}). Let us write the boundary divisor as $D+\E$, where $\E$ is the reduced exceptional divisor of the resolution of $S'$. Using the fact that $\ovk(S_0)=0$ and that the intersection matrix of $D$ is not negative definite we get some restrictions on the shape of $D$ following from \cite[8.8]{Fujita}. In fact for smooth exceptional surfaces this would be enough to get the description of them. However, in the singular case we need to obtain more restrictions on $D$ because we do not have much control over $\E$. We do this in section \ref{s:rulings}. In remaining ten cases we are able to find $0$-curves inside $D$, which give $\PP^1$-rulings of the completion having nice properties. We analyze singular fibers and sections of these rulings and we eliminate all but two cases. Having enough information on the latter two rulings in section \ref{s:constructions} we are able to construct two exceptional singular $\Q$-homology planes and prove their uniqueness. We compute their automorphism groups, the orders of the first homology groups and show that they came from special line arrangements on $\PP^2$.

\tableofcontents

\section{Preliminaries}

For convenience we recall some facts from the theory of open algebraic surfaces that we use more often, partially to fix the notation. The reader is referred to \cite{Miyan-OpenSurf} for details.

\subsection{\label{ss:generalities}Divisors}

Let $D=\sum_{i=1}^n m_iD_i$ with $D_i$ distinct, irreducible and $m_i\in\Q\setminus\{0\}$ be a simple normal crossing (snc-) divisor on a smooth complete surface. Put $$d(D)=det(-Q(D)),$$ where $Q(D)$ is the intersection matrix of $D$, i.e. $Q(D)_{i,j}=m_im_jD_i\cdot D_j$. We define the reduction of $D$ as $\un D=\sum D_i$ and denote the number of components of $D$ by $\#D$. By a component we always mean an irreducible component. The numerical equivalence of divisors is denoted by $\equiv$. We write $D\geq 0$ for effective divisors and for $\Z$-divisors linearly equivalent to effective divisors. Two $\Q$-divisors $A,B$ are linearly equivalent if $rA$ and $rB$ are linearly equivalent $\Z$-divisors for some nonzero integer $r$. For a $\Q$-divisor $D$ linearly equivalent to some effective $\Q$-divisor we write $D\geq_{\Q}0$.

The dual graph of $D$ is a weighted one-dimensional simplicial complex with one vertex $v_i$ of weight $D_i^2$ for each irreducible component $D_i$ of $D$ and one edge between $v_i$ and $v_j$ for each point of intersection of $D_i$ with $D_j$. We say that $D$ is a \emph{forest} (\emph{tree}) if $\Supp D$ is simply connected (and connected). It is \emph{rational} if all its components are rational. $D$ is a \emph{chain} if it is connected and each component $D_i$ of $D$ is non-branching, i.e. it has the \emph{branching number} $\beta_D(D_i)=D_i\cdot (D-D_i)$ not greater than two. A \emph{tip} is a component with $\beta_D\leq 1$. A chain $D$ is \emph{admissible} if it is rational and $D_i^2\leq-2$ for every $i$. A curve $L$ is a \emph{$(b)$-curve} if and only if $L\cong\PP^1$ and $L^2=b$. We say that $D$ is a \emph{fork} if it is a tree with a unique branching component $B$ and $\beta_D(B)=3$. Suppose $R$ is a rational chain with some tip $R_1$ chosen. We write $$R=[-R_1^2,-R_2^2,\ldots,-R_r^2],$$ where $R_i$'s are components of $R$ ordered in such a way that $R_i\cdot R_{i+1}=1$ for $i=1,\ldots,r-1$ and we define $d'(R)=d(R-R_1)$ with $d(0):=1$. $R^t$ is the same chain as $R$ but considered with a reversed order. If $R$ is a $(-2)$-chain, i.e. $R=[2,2,\ldots,2]$, then we write $R=[(r)]$, where $r=\#R$. If $R$ is admissible we define $$\delta(R)=\frac{1}{d(R)},\ e(R)=\frac{d'(R)}{d(R)},\ \wt e(R)=e(R^t).$$ If $D$ is not a chain we define its \emph{maximal twigs} as the rational chains of maximal length with support contained in $\Supp D$, which do not contain branching components of $D$ and contain a tip of $D$. Each twig is considered with a natural linear order on the set of components for which its tip is the first component. If $D$ is not an admissible chain we define its \emph{maximal admissible twigs}, say $T_1,\ldots,T_s$, analogously and put $$\delta(D)=\sum_{i=1}^s\delta(T_i),\ e(D)=\sum_{i=1}^s e(T_i),\ \wt e(D)=\sum_{i=1}^s\wt e(T_i).$$

\emph{Smooth pair} $(X,D)$ consists of a smooth complete (hence projective by the result of Zariski) surface and a reduced snc-divisor on it. In this case we write $X-D$ for $X\setminus \Supp D$. The divisor $D$ is \emph{snc-minimal} if after a contraction of any $(-1)$-curve in $D$ the direct image of $D$ is not an snc-divisor. A smooth pair $(X,D)$ is snc-minimal if $D$ is snc-minimal. The pair $(X,D)$ is a \emph{smooth completion} of an open surface $U$ if $X-D=U$.

If $\pi:X'\to X$ is a birational morphism then we write $\pi^{-1}(D)$ for the \emph{preimage} of $D$, which we define as $\un {\pi^*D}$, the reduced total transform of $D$. A blowup with a center on an snc-divisor $D$ is \emph{subdivisional} for $D$ if the center belongs to two components of $D$, otherwise it is \emph{sprouting for $D$}. The sequence of blowups over $D$ (i.e. with centers on $D$ and on its successive preimages) is subdivisional if all blowups are subdivisional for the respective preimages of $D$. The composable sequence of blowups is \emph{connected} if the exceptional divisor of the composition contains a unique $(-1)$-curve.

\subsection{\label{ss:rulings}Rulings}

We say that a surface $X$ is \emph{$\PP^1$-ruled} (respectively \emph{$\C^1$-ruled}, \emph{$\C^*$-ruled}, \emph{$\C^{n*}$-ruled}) if there exists a curve $B$ and a surjective morphism $p:X\to B$ with a general fiber isomorphic to $\PP^1$ (respectively to $\C^1$, $\C$ with one, $\C$ with $n+1$ points deleted). We call also the $\C^1$-ruling an \emph{affine} ruling. Clearly, if $X$ is normal then $B$ can be assumed to be smooth.

Suppose that $X$ is smooth and has a ruling as above. Then for some smooth completion $(\ov X,D)$ this ruling can be extended to a $\PP^1$-ruling $\overline{p}:\ov X\to\overline{B}$, where $\overline{B}$ is a smooth completion of $B$. Let $F$ be a fiber of $\overline{p}$. An irreducible curve $C\subseteq \ov X$ is called an \emph{$n$-section} if $F\cdot C=n$. We will say just \emph{section} for a $1$-section. $C$ is \emph{horizontal} if $n>0$, otherwise it is \emph{vertical}. If $C$ is vertical then it is called a \emph{$D$-component} if $C\subseteq D$, otherwise it is called an \emph{$X$-component}. If the ruling is fixed we denote the divisor consisting of horizontal components of $D$ by $D_h$. The divisor is horizontal (vertical) if all its components are horizontal (vertical). The completion $(\ov X,D)$ is \emph{$p$-minimal} if it is smooth and minimal with respect to the property that the extension of $p$ from $X$ to $\ov X$ exists (the partial order is induced by morphisms of pairs).

For a smooth pair $(\ov X,D)$ put $X=\ov X-D$. Let $\pi$ be a $\PP^1$-ruling of $\ov X$. Following \cite{Fujita} we define some characteristic numbers of the triple $\tau=(\ov X,D,\pi)$: $h_\tau$ is the number of horizontal $D$-components, $\sigma_\tau(F)$ is the number of $X$-components contained in $F$, $\nu_\tau$ is the number of fibers contained in $D$. Put $$\Sigma_\tau=\underset{F\nsubseteq D}{\sum}(\sigma_\tau(F)-1).$$ If there is no danger of confusion we omit indices writing $\Sigma$ (or $\Sigma_X$) for $\Sigma_\tau$, $h$ for $h_\tau$, etc. If one contracts a vertical $(-1)$-curve and simultaneously changes $\ov X$ and $D$ for their images then the numbers $b_2(\ov X)-b_2(D)-\Sigma+\nu$ and $h$ do not change ($b_i(X)=\dim H_i(X;\Q)$). This leads to the following equation (cf. \cite[4.16]{Fujita}):

\begin{align} \Sigma=h+\nu+b_2(\ov X)-b_2(D)-2. \label{eq:Fujita}\end{align}

\noin Clearly, $b_2(\ov X)-b_2(D)$ depends only on $X$.

We now summarize some information about singular fibers of $\PP^1$-rulings (cf. \cite[\S 4]{Fujita}). For a given ruling $\pi$ and a vertical component $C$ the multiplicity $\mu(C)$ is the coefficient of $C$ in $\pi^*(\pi(C))$.

\blem\label{lem:singular fibers} Let $F$ be a singular fiber of a $\PP^1$-ruling of a smooth complete surface. Then $F$ is a rational snc-tree containing a $(-1)$-curve. Each $(-1)$-curve of $F$ intersects at most two other components of $F$. Successive contractions of $(-1)$-curves contract $F$ to a smooth $0$-curve. In this process the number of $(-1)$-curves can increase only in the last but one step, when $[2,1,2]$ contracts to $[1,1]$.\elem

Suppose that $F$ as above contains a unique $(-1)$-curve $C$. The sequence of blowups recovering $F$ from a smooth $(0)$-curve is connected. Let $B_1,\ldots, B_n$ be the branching components of $F$ written in order in which they are produced in the sequence of blowups recovering $F$ from a smooth $(0)$-curve and let $B_{n+1}=C$. We can write $\un F$ as $\un F=T_1+T_2+ \ldots+T_{n+1}$, where the divisors $T_i$ are chains consisting of all components of $\un F-T_1-\ldots -T_{i-1}$ created not later than $B_i$. We call $T_i$ the \emph{i-th branch} of $F$ and say that $F$ is \emph{branched} if $i>1$.

\brem\label{rem:singular fibers} Let $F$ and $C$ be as above. Then $\mu(C)>1$ and there are exactly two components of $F$ having multiplicity one. They are tips of the fiber and belong to the first branch. The connected component of $\un F-C$ not containing curves of multiplicity one is a chain. If $\mu(C)=2$ then either $F=[2,1,2]$ or $C$ is a tip of $F$ and then $\un F-C$ is either a $(-2)$-chain or a $(-2)$-fork with two tips as maximal twigs.\erem

\subsection{\label{ss:barks}Zariski decomposition}

Let $(X,D)$ be a smooth pair. If it is almost minimal (cf. \cite[2.3.11]{Miyan-OpenSurf}) and $\ovk(X-D)\geq 0$ then the Zariski decomposition of $K_X+D$, where $K_X$ stands for the canonical divisor on $X$, can be computed explicitly using a \emph{bark of $D$}. For non-connected $D$ bark is a sum of barks of its connected components, so we will assume $D$ is connected. If $D$ is an snc-minimal resolution of a quotient singularity (i.e. $D$ is an admissible chain or an admissible fork, cf. \cite[2.3.4]{Miyan-OpenSurf}) then we define $\Bk D$ as a unique $\Q$-divisor with $\Supp \Bk D\subseteq D$, such that  $$(K_X+D-\Bk D)\cdot D_i=0\text{\ for each component\ } D_i\subseteq D.$$ In other case let $T_1,\ldots, T_s$ be all the maximal admissible twigs of $D$. (If $\ovk(X-D)\geq 0$ and $D$ is snc-minimal then all maximal twigs of $D$ are admissible, cf. \cite[6.13]{Fujita}). We define $\Bk D$ as a unique $\Q$-divisor with $\Supp \Bk D\subseteq \bigcup T_j$, such that $$(K_X+D-\Bk D)\cdot D_i=0\text{\ for each component\ } D_i\subseteq \bigcup_{j=1}^s T_j.$$

Suppose $R$ is an admissible chain with some tip $R_1$ chosen. Then we define $\Bk (R,R_1)$ as a unique $\Q$-divisor with support contained in $R$, such that $$R_1\cdot \Bk (R,R_1)=-1\text{ and }\ R_i\cdot \Bk (R,R_1)=0 \text{\ for\ each component\ } R_i\subseteq R-R_1.$$ If there is no need to mention the tip explicitly (for example if $R$ is an admissible twig of some fixed divisor then its tip will be a default choice for $R_1$) we write $\Bk' R$ instead of $\Bk (R,R_1)$.

(This notation does not occur in standard references, but we find it useful). Now we can write $\Bk D=\Bk' T_1+\ldots+\Bk' T_s$. We recall here the properties of $\Bk D$ which we use later and refer the reader to \cite[\S 2.3]{Miyan-OpenSurf} for details. We put $D^\#=D-\Bk D$.

\blem\label{lem:Bk properties} Let $(X,D)$ be a smooth pair. Write $D=\sum D_i$ with $D_i$ distinct irreducible and $\Bk D=\sum d_i D_i$. One has:\benum[(i)]

\item $0\leq d_i\leq 1$ for each $i$, $\Bk D$ is rational and $Q(\Bk D)$ is negative definite, unless $\Bk D=0$,

\item if $d_i=1$ for some $i$ and $D'$ is a connected component of $D$ containing $D_i$ then $\Bk D'=D'$ and $D'$ consists of $(-2)$-curves,

\item $\Supp \Bk D$ consists of the supports of all maximal admissible twigs of $D$ and of all connected components of $D$ which are either admissible chains or admissible forks (see \cite[2.3.5]{Miyan-OpenSurf}),

\item $(K_X+D^\#)\cdot  Z=0$ for every $Z\subseteq \Supp \Bk D$,

\item if $(X,D)$ is almost minimal and $\ovk(X-D)\geq 0$ then $(K_X+D)^-=\Bk D$.

\eenum\elem

We now state a version of Bogomolov-Miyaoka-Yau inequality proved by Langer (\cite[Corollary 5.2]{Langer}), which generalizes the inequalities of Miyaoka \cite[Theorem 1.1]{Miyaoka} and Kobayashi \cite[Theorem 2]{Kob}. See \cite[3.4, \S9]{Langer} for a definition of the \emph{orbifold Euler number} $\chi_{orb}(X,D)$ and for computations in special cases.

\bprop\label{thm:Langer-ineqality} Let $(X,D)$ be a normal projective surface together with a $\Q$-divisor $D=\sum m_iD_i$ with $0\leq m_i\leq 1$. Assume that the pair is log-canonical and $K_X+D$ is pseudoeffective. Then $$3 \chi_{orb}(X,D)+\frac{1}{4}((K_X+D)^-)^2\geq (K_X+D)^2.$$ \eprop

\bcor\label{lem:KobIneq} Let $(X,D)$ be a smooth pair with $\kp(K_X+D)\geq 0$. Then:\benum[(i)]

\item $$3\chi(X-D)+\frac{1}{4}((K_X+D)^-)^2\geq (K_X+D)^2.$$

\item For each connected component of $D$, which is a connected component of $\Bk D$ (hence contractible to a quotient singularity) denote by $G_P$ the local fundamental group of the respective singular point $P$. Then $$\chi(X-D)+\sum_P\frac{1}{|G_P|}\geq \frac{1}{3} (K_X+D^\#)^2.$$ \eenum \ecor

\begin{proof} According to \cite[7.6]{Langer} if $(X,D)$ is a pair as in \ref{thm:Langer-ineqality} and $D$ is reduced then for a point $P\in D$ the local orbifold numbers $\chi_{orb}(P;X,D)$ vanish, hence $$\chi_{orb}(X,D)=\chi(X-\Sing X-D)+\sum_{P\in\Sing X}\chi_{orb}(P;X,D).$$ This already proves (i), where $X$ is smooth. Let $\pi:(X,D)\to (X',D')$ be a morphism contracting the connected components of $\Bk D$ to quotient singularities. Then by \cite[2.3.14.1]{Miyan-OpenSurf} $K_X+D^\#\equiv \pi^*(K_{X'}+D')$ and $K_{X'}+D'=\pi_*(K_X+D^\#)$, in particular $K_{X'}+D'$ is pseudoeffective because $(K_X+D)^--\Bk D$ is effective by \ref{lem:Bk properties}(iv) and the properties of the Zariski decomposition of $K_X+D$. We need to know $\chi_{orb}(P;X',D')$. If $P\not\in D'$ then the preimage of $P$ is a connected component of $D$ (and of $\Bk D$) and by \cite[7.1]{Langer} we have $\chi_{orb}(P;X',D')=\frac{1}{|G_P|}$. We have also $\chi(X'-\Sing X'-D')=\chi(X-D)$. Since $((K_X'+D')^-)^2\leq 0$, (ii) follows from \ref{thm:Langer-ineqality} applied to $(X',D')$.\end{proof}

\bsrem  Part (ii) generalizes the Kobayashi inequality for the case $\ovk(X-D)=0,1$, it is stronger than the original Miyaoka inequality (there is no $\frac{1}{4}N^2$ term, using the notation of \cite[Theorem 1.1]{Miyaoka}). If $\ovk(X-D)=2$ then to get the original Kobayashi inequality one applies \ref{thm:Langer-ineqality} to the \emph{strongly minimal model} of $(X,D)$ (cf. \cite[2.4.12, 2.6.6]{Miyan-OpenSurf}). \esrem

\subsection{Other useful results}

As a consequence of elementary properties of determinants one gets the following result.

\blem\label{lem:DetFormulas}(\cite[2.1.1]{KR-ContrSurf}). Let $D$ be a reduced snc-tree. \benum[(i)]

\item Let $C$ be a component of $D$ and let $D_1,D_2,\ldots,D_k$ be the connected components of $D-C$. If $C_i$ is the component of $D_i$ meeting $C$ then $$d(D)=-C^2\prod_i d(D_i)-\sum_i d(D_i-C_i)\prod_{i\neq j}d(D_j).$$

\item Let $D=D_1+D_2$, where $D_1, D_2$ are connected and intersect in one point. Let $C_1\subseteq D_1,C_2\subseteq D_2$ be the intersecting components, then $$d(D)=d(D_1)d(D_2)-d(D_1-C_1)d(D_2-C_2).$$ \eenum\elem

\bsrem If $D$ is an snc-divisor then $d(D)$ is invariant under blowup, i.e. if $(X,D)$ is a smooth pair and $\sigma:X'\to X$ is a blowup, then $d(\sigma^{-1}(D))=d(D)$. For trees this follows from \ref{lem:DetFormulas} by induction on $\#D$. \esrem

\blem\label{lem:Eff-NegDef}
Let $A$ and $B$ be some $\Q$-divisors, such that $A+B$ is effective and $Q(B)$ is negative definite. If $A\cdot B_i=0$ for each irreducible component $B_i$ of $B$ then $A$ is effective.\elem

\begin{proof} We can assume that $A$ and $B$ are $\Z$-divisors and $B$ is effective and nonzero. Write $B=\sum b_iB_i$ for some positive integers $b_i$ and irreducible components $B_i$ of $B$. Choose $b_i'\in \N$, such that the sum $\sum b_i'$ is the smallest possible among divisors $\sum b_i'B_i$, such that $A+\sum b_i'B_i$ is effective. If $b_i'>0$ for some $i$ then $(A+\sum b_i'B_i)\cdot (\sum b_i'B_i)=(\sum b_i'B_i)^2<0$ by the assumptions. Hence $\Supp (A+\sum b_i'B_i)$ contains some $B_i$, a contradiction with the definition of $b_i'$. Thus $A$ is effective. \end{proof}

\blem\label{lem:chi of minimal model is smaller} Let $X_0$ be a smooth part of $X'\setminus D'$, where $D'$ is a divisor on an affine surface $X'$. Let $(X_m,D_m)$ be the almost minimal model of some smooth completion of $X_0$. Then the almost minimal model $X_m-D_m$ of $X_0$ is an open subset of $X_0$ and $\chi(X_m-D_m)\leq \chi(X_0)$.\elem

\begin{proof} Let $\epsilon:X\to X'$ be a resolution with snc-minimal exceptional divisor and let $(\ov X,D)$ be a smooth completion of $X$. Since $X'$ is affine, $D$ is connected and $Q(D)$ is not negative definite. Let $D''\subseteq \ov X$ be the closure of $\epsilon^{-1}(D')$ and let $E$ be the part of the exceptional divisor with support equal to $\epsilon^{-1}(\Sing (X'-D'))$. Blowing on $D''$ if necessary we can assume that $(\ov X,D+D''+E)$ is a smooth completion of $X_0$. Moreover, $D+D''$ is connected. Consider the process of producing an almost minimal model $(X_m,D_m)$ of $(\ov X,D+D''+E)$, it goes by contractions of special $(-1)$-curves, so-called \emph{log-exceptional curves of the first kind} (cf. \cite[2.4.3]{Miyan-OpenSurf}). Notice that in the process the divisor $D'+D''$ cannot be contracted, because $Q(D)$ is not negative definite. By the properties of a log-exceptional curve not contained in the boundary its contraction causes a subtraction of a curve with $\chi=1$ or $\chi=0$ from $X_0$. Contractions of $(-1)$-curves contained in the boundary divisor do not affect $X_0$, unless some connected component of the boundary is eventually contracted to a smooth point which does not belong to the proper image of the boundary divisor. Then this point adds to the almost minimal model of $X_0$. Affiness of $X'$ implies that a log-exceptional curve not contained in $E$ intersects the image of $D$, so the above cannot happen for connected components of $E$.\end{proof}

\section{Basic properties of $S'$}\label{s:basic_properties}

We now fix the notation for the rest of the paper. Let $S'$ be an exceptional singular $\Q$-homology plane, i.e. its smooth locus $S_0$ has $\ovk(S_0)\neq 2$ and is neither $\C^1$- nor $\C^*$-ruled. As was explained in section \ref{s:main result}, this implies $\ovk(S_0)=0$. Let $\epsilon \colon S \to S'$ be a resolution having an snc-divisor as the exceptional locus and let $(\ov S,D)$ be a smooth completion of $S$. By the definition of the logarithmic Kodaira dimension $\ovk(S')=\ovk(S)=\kp(K_{\ov S}+D)$, where $K_{\ov S}$ stands for the canonical divisor on $\ov S$. Let $\{p_1, \ldots, p_q\}$ be the singular locus of $S'$ and let $\E_i=\epsilon^{-1}(p_i)$. We assume that $\E =\E_1+\E_2+\ldots +\E_q$ is snc-minimal. The intersection matrix $Q(\E)$ is negative definite. We write $H_i(X,A)$ for $H_i(X,A;\Q)$ and $b_i(X,A)$ for $\dim H_i(X,A;\Q)$.

\blem\label{lem:topology} Let $i:D\cup\E \to \ov S$ be the inclusion. The following properties hold:\benum[(i)]

\item $H_2(i):H_2(D\cup\E)\to H_2(\ov S)$ is an isomorphism,

\item $S'$ is rational,

\item $D$ is a rational tree,

\item $\Sigma_{S_0}=h+\nu-2$ and $\nu\leq 1$,

\item $S'$ is affine.\eenum\elem

\begin{proof} (i) Let $Tub(\E)$ be a sum of tubular neighborhoods of $\E_i$'s in $S$ (see \cite{Mumford} for the construction) and let $M$ be the boundary of the closure of $Tub(\E)$. We can assume that $M$ is a disjoint sum of closed oriented 3-manifolds. There exists a deformation retraction $Tub(\E)\to \E$, so by excision $H^j(S_0,M)=H^j(S,Tub(\E))=H^j(S,\E)$ and since for $j>1$ we have $H^j(S,\E)=H^j(S')=0$, we get $b_j(S_0)=b_j(M)$ for $j>1$. In fact $b_1(S_0)$ also equals $b_1(M)$ because $H^1(S_0,M)=H^1(S,\E)=\Q^{q-1}$ and then $H^0(M)\to H^1(S_0,M)$ is an epimorphism. By \cite{Mumford} $b_1(M)=b_1(\E)=0$, so each connected component of $M$ is a $\Q$-homology sphere by the Poincare duality. We conclude that $b_j(S_0)=0$ for $j=1,2$. Now by the Lefschetz duality $H_j(\ov S,D\cup\E)=H^{4-j}(S_0)$, hence $H_2(\ov S,D\cup\E)=H_3(\ov S,D\cup\E)=0$. It follows from the exact sequence of the pair $(\ov S,D\cup\E)$ that $H_2(i)$ is an isomorphism.

(ii) Since $H_2(i)$ is an isomorphism, the exact sequence of the pair $(\ov S,D\cup\E)$ gives that $H_3(\ov S)\to H_3(\ov S,D\cup\E)$ is an isomorphism. Therefore by the Lefschetz duality $b_1(\ov S)=b_3(\ov S)=b_3(\ov S,D\cup\E)=b_1(S_0)=0$. Now if $\kp(\ov S)=-\8$ then $\ov S$ is birational to a $\PP^1$-fibration over some complete curve $B$. From the homotopy exact sequence of a fibration we know that $b_1(B)=b_1(\ov S)$, so $B\cong \PP^1$, hence $\ov S$ is rational. Suppose $\kp(\ov S)\geq 0$. Since $\kp(\ov S)\leq \ovk(S_0)=0$, we see that $\kp(\ov S)=\ovk(S)=\ovk(S_0)=0$. We now prove that $Q(D)$ is negative definite. We can assume that $(\ov S,D)$ is almost minimal. Then by \ref{lem:Bk properties}(v) $K_{\ov S}+D^\#=(K_{\ov S}+D)^+\equiv 0$ and $K_{\ov S}\geq_{\Q} 0$, so $D^\#=0$ because $D^\#$ is effective. Thus $D=\Bk D$, and we are done by \ref{lem:Bk properties}(i). By (i) we get a contradiction with the Hodge index theorem.

(iii) Since $H_2(\ov S,D\cup\E)=0$, the exact sequence of the pair $(\ov S,D\cup\E)$ gives the injectivity of $H_1(D\cup\E)\to H_1(\ov S)$, so $b_1(D)=0$ by (ii). In the proof of (i) we have shown that $b_1(\ov S,D\cup\E)=b_3(M)$, so since $M$ is a disjoint sum of $b_0(\E)$ three-dimensional manifolds, we get $b_1(\ov S,D\cup\E)=b_0(\E)$. Now the exact sequence of a pair $(\ov S,d\cup \E)$ gives $b_0(D\cup\E)=b_1(\ov S,D\cup\E)+b_0(\ov S)=b_0(\E)+1$, hence $D$ is connected.

(iv) The first equation is a consequence of \eqref{eq:Fujita} and (i). If $\nu>1$ then the numerical equivalence of fibers of a $\PP^1$-ruling gives a numerical dependence of components of $D+\E$ in $NS(\ov S)\otimes\Q$, where $NS(\ov S)$ is the Neron-Severi group of $\ov S$. This contradicts (i).

(v) Since $H_2(i)$ is an epimorphism by (i) and since $D$ is connected by (iii), Fujita's argument from the proof of \cite[2.4(3)]{Fujita} works.

\end{proof}

\bsrem From \ref{lem:topology}(i) and the Hodge index theorem we get $d(D+\E)<0$, so $d(D)<0$.\esrem

\blem\label{lem:S0-L generic} Every irreducible curve $L\nsubseteq D\cup\E$ satisfies $\ovk(S_0-L)=2$. \elem

\begin{proof}
Suppose $\ovk(S_0-L)=1$. Since $S_0$ does not contain complete curves, \cite[2.3]{Kawamata2} implies that $S_0-L$ is $\C^*$-ruled. $S_0$ is not $\C^*$-ruled, so it is affine-ruled and we get $\ovk(S_0)=-\8$ by the easy addition theorem (\cite[Theorem 10.4]{Iitaka}), a contradiction. Suppose $\ovk(S_0-L)=0$. By \ref{lem:topology}(i) $H_2(\ov S,\Q)$ is generated by cycles contained in $D\cup\E$, hence $NS(\ov S)\otimes\Q$ is generated by the components of $D+\E$. Since $\ov S$ is rational, we get $Pic(S_0)\otimes\Q=0$, so there exists a rational function $f$ on $S_0$, such that $(f)=kL$ for some $k>0$. We get a morphism $f:S_0-L\to \C^*$. If $S_0-L\to B\to\C^*$ is its Stein factorization then $\ovk(B)\geq \ovk(\C^*)=0$ and $0\geq \ovk(F_b)+\ovk(B)$ for a fiber $F_b$ over a generic $b\in B$ by Kawamata addition theorem (\cite{Kawamata}). Since $S_0-L$ is not affine ruled, we get $\ovk(F_b)=0$, i.e. $f$ is a $\C^*$-ruling, a contradiction.
\end{proof}

\bdfn Let $(X,B)$ be a smooth pair. A curve $C\subseteq X$ is a \emph{simple curve on $(X,B)$} if and only if $C\cong \PP^1$ and $C$ has at most one common point with each connected component of $B$.\edfn

\bcor\label{cor:no simple curve} There is no simple curve on $(\ov S,D+\E)$. If $D$ is snc-minimal then the pair $(\ov S,D+\E)$ is almost minimal.\ecor

\begin{proof}
Let $L$ be a simple curve on $(\ov S,D+\E)$. Since $S'$ is affine, $L\cap D\neq\emptyset$. Let $(X_m,B_m)$ be the almost minimal model of some smooth completion of $S_0-L$ and let $(X_m,B_m)\to (X_r,B_r)$ be the morphism contracting the connected components of $\Bk B_m$. Denote the local fundamental group of a singular point $P\in\Sing(X_r-B_r)$ by $G_P$. By \ref{lem:chi of minimal model is smaller} $X_m-B_m$ is an open subset of $S_0-L$ satisfying $\chi(X_m-B_m)\leq \chi(S_0-L)$. Since $(X_m,B_m)$ is almost minimal, by \ref{lem:Bk properties}(v) $(K_{X_m}+B_m)^+\equiv K_{X_m}+B_m^\#$, so by \ref{lem:KobIneq}(ii) and \ref{lem:S0-L generic} $\chi(X_m-B_m)+\sum\frac{1}{|G_P|}>0$. Put $s=|L\cap\E|$. The matrix $Q(D)$ is not negative definite, so $|\Sing(X_r-B_r)|\leq q-s$. This gives $\sum\frac{1}{|G_P|}\leq \frac{q-s}{2}$, so $\chi(S_0-L)\geq\chi(X_m-B_m)> -\sum\frac{1}{|G_P|}\geq\frac{s-q}{2}$. We compute $\chi(S_0-L)=\chi(S_0)-\chi(L)+|L\cap D|+s=1-q+s-2+|L\cap D|$, hence $|L\cap D|=\chi(S_0-L)+1+q-s>\frac{q-s}{2}+1$, so $|L\cap D|>1$, a contradiction. Since log-exceptional curves of the first kind not contained in $D\cup \E$ are simple, $(\ov S,D+\E)$ is almost minimal. \end{proof}

By \ref{lem:topology}(iii) $D$ is a rational tree and since $Q(D)$ is not negative definite, if it is snc-minimal then by \cite[8.8]{Fujita} it is of one of the following types:\benum[]

\item (Y): a fork with three maximal admissible twigs and $\delta(D)=1$,

\item (H): has dual graph $$ \xymatrix{{-2}\ar@{-}[r] &{\cdot}\ar@{-}[r]\ar@{-}[d]& {\cdots}\ar@{-}[r] &{\cdot}\ar@{-}[r]\ar@{-}[d]&{-2}\\ {} & {-2} & {} & {-2} & {} } ,$$

\item (X): has dual graph $$\xymatrix{{} & {-2}\ar@{-}[d] & {} \\ {-2}\ar@{-}[r] &{\cdot}\ar@{-}[r]\ar@{-}[d]& {-2} \\
{} & {-2} & {} }$$\eenum

\ssk We will frequently use the fact that, as a consequence of the Riemann-Roch theorem, on a complete rational surface a $(0)$-curve (and hence any rational tree which contracts to a $(0)$-curve) induces a $\PP^1$-ruling with this curve as one of the fibers.

\section{Rulings of $S_0$ with $\nu>0$}\label{s:rulings}

From now on we assume that $D$ is snc-minimal.

\blem\label{lem:S_0-components are exceptional} Let $D_0\subseteq D$ be a component of $D$ meeting some maximal twig of $D$ and such that $D_0^2\geq 0$. Let $\sigma:(\wt S,\wt D)\to (\ov S, D)$ be modification over $D$ obtained by blowing up successively in the point of intersection of $D_0$ with the preimage of this maximal twig until $D_0^2=0$. Let $\pi:\wt S\to \PP^1$ be the induced $\PP^1$-ruling with $D_0$ as a fiber. Then a component of a fiber is an $S_0$-component if and only if it is exceptional.\elem

\begin{proof}Denote the maximal twig of $D$ as above by $T$. Let $L$ be an $S_0$-component of some fiber. We have $\ovk(S_0)=0$, so $(K_{\ov S}+D+\E)^+\equiv 0$ by \cite[6.11]{Fujita} and then $K_{\ov S}+D+\E\equiv \Bk D+\Bk \E$ by \ref{cor:no simple curve} and \ref{lem:Bk properties}(v). The sequence of blowups defining $\sigma$ is subdivisional for $D$, so $K_{\wt S}+\wt D+\E\equiv \sigma^*\Bk D+\Bk \E$ and $L^2=-2-L\cdot K_{\wt S}=-2+L\cdot (\wt D-\sigma^*\Bk D)+L\cdot (\E-\Bk \E)\geq -2+L\cdot (\wt D-\sigma^*\Bk D)$. Since $D_0\nsubseteq \Supp \Bk D$ by \ref{lem:Bk properties}(i) and each of the blowups is sprouting for the respective preimages of $T$, by \ref{lem:Bk properties}(ii) the coefficients of components of $\wt D$ in $\sigma^*\Bk D$ are smaller than one. Thus $L^2>-2$ because $L\cdot \wt D>0$, so we are done.\end{proof}

\bsrem Notice that no fiber of a $\PP^1$-ruling of $\wt S$ can be contained in $\E$, otherwise $\wt D$ would be vertical, so $S'$ would contain complete curves.\esrem

The following lemma, which is a generalization of arguments from \cite[6.2]{Koras-A2} allows to bound from below the self-intersection of one of the branching components of $D$ having four maximal twigs.

\blem\label{lem:table} Let $T$ be an snc-minimal divisor with two branching components $B$, $B'$ and such that $\beta_T(B)=\beta_T(B')=3$. Let $T_1$, $T_2$ and $T_3$, $T_4$ be the maximal twigs of $T$ intersecting $B$ and $B'$ respectively. If $Q(T-B-B')$ is negative definite, $\wt e(T_1)+\wt e(T_2)\leq -B^2-1$ and $\wt e(T_3)+\wt e(T_4)\leq -B'^2-1$ then either $Q(T)$ is negative definite or $d(T)=0$ and then $T-T_1-T_2-T_3-T_4$ is a $(-2)$-chain and $\wt e(T_1)+\wt e(T_2)=\wt e(T_3)+\wt e(T_4)=1$.\elem

\begin{proof} Write $T-T_1-T_2-T_3-T_4=B_1+B_2+\ldots+B_t$ with $B_1=B$ and $B_t=B'$. Define $T_0=B_2+\ldots+B_{t-1}$ and $d_i=d(D^{(i)})$, where $D^{(i)}=T_3+T_4+B_t+B_{t-1}+\ldots+B_i$. By \ref{lem:DetFormulas}(i) $d_2=d(T_3)d(T_4)d(T_0)(-B_t^2-\wt e(T_3)-\wt e(T_4)-\wt e(T_0))\geq d(T_3)d(T_4)d(T_0)(1-\wt e(T_0))>0$, so $D^{(2)}$ is negative definite by Sylvester's theorem. We now prove that $d_2\geq d_3$. Notice that $B_i^2\leq -2$ for each $i$, so by \ref{lem:DetFormulas} for $i=2,\ldots,t-1$ we get $d_{i}-d_{i+1}=(-B_{i}^2-2)d_{i+1}+d_{i+1}-d_{i+2}\geq d_{i+1}-d_{i+2}$. We have $d_{t}-d_{t+1}=d(T_3)d(T_4)(-B_t^2-\wt e(T_3)-\wt e(T_4)-1)\geq 0$, so we are done. By \ref{lem:DetFormulas}(i) $d(T)=d_2d(T_1)d(T_2)(-B_1^2-\frac{d_3}{d_2}-\wt e(T_1)-\wt e(T_2))\geq d_2d(T_1)d(T_2)(-B_1^2-\wt e(T_1)-\wt e(T_2)-1)\geq 0$. Hence by Sylvester's theorem if $d(T)\neq 0$ then $Q(T)$ is negative definite. On the other hand, if $d(T)=0$ then all the weak inequalities above become equalities and the thesis follows.\end{proof}

From now on we denote the maximal twigs of $D$ by $T_1,\ldots,T_s$. If $D$ has only one branching component we denote it by $B$.

\blem\label{lem:boundary description} $D$ can be only of type (X) or (Y). If it is of type (X) then $-1\leq B^2\leq 0$. If it is of type (Y) then $B^2=-1$ and the triple $(d(T_1),d(T_2),d(T_3))$ is up to permutation one of the following: $(3,3,3)$, $(2,3,6)$, $(2,4,4)$.\elem

\begin{proof} In case (H) let $B, B'$ be the branching components of $D$. The chain $D-T_1-T_2-T_3-T_4-B-B'$ is admissible, otherwise after some subdivisional modification of $D$ it contains a $(0)$-curve, hence gives a $\C^*$-ruling of $S_0$, which contradicts our assumptions about $S_0$. Since $Q(D)$ is not negative definite and $d(D)\neq 0$, by \ref{lem:table} we can assume that $B^2\geq -1$. Assume $T_1$ and $T_2$ meet $B$. Blow up on the intersection of $B$ with $D-T_1-T_2-B$ until $B^2=-1$. We have $T_1^2=T_2^2=-2$, so $T_1+2B+T_2$ gives a $\C^*$-ruling of $S_0$, a contradiction. Thus only types (X) and (Y) remain. We have $d(D)<0$, so by \ref{lem:DetFormulas} $-B^2+\delta(D)<\sum_{i=1}^s \frac{d'(T_i^t)+1}{d(T_i)}\leq s$. For both types we obtain $B^2\geq -1$.

In case (Y) we have $\delta(D)=1$ by definition, so we need only to prove that $B^2=-1$. Suppose $B^2>0$ in case (X) or $B^2\geq 0$ in case (Y). Let $\sigma:(\wt S,\wt D)\to (\ov S,D)$ be the modification obtained by blowing up the point of intersection of $T_1$ with $B$ until $B^2=0$. Consider the $\PP^1$-ruling of $\wt S$ given by $B$. We see that $\wt D$ contains no vertical $(-1)$-curves. The divisor $D_h$ consists of three or four sections of the ruling. Put $D_v=\wt D-D_h-B$. Notice that if some section intersects a vertical component $V$ then $\mu(V)=1$ and it does not intersect any other component lying in the fiber containing $V$. By \ref{lem:S_0-components are exceptional} the $S_0$-components of singular fibers are exceptional.

Let $F$ be a fiber containing some connected component of $\E$ ($\E$ is vertical, because $B\cdot \E=0$). If $F$ contains some $\wt D$-components then there exists a chain of $S_0$-components in $F$ connecting $\E \cap F$ with some $\wt D$-component of $F$. In fact this chain consists of a unique $(-1)$-curve $L$, since all $S_0$-components are $(-1)$-curves and two of them cannot meet. By \ref{cor:no simple curve} $D_h\cdot L>0$, so $\mu(L)=1$, a contradiction. Therefore there are no $\wt D$-components in $F$, hence each $S_0$-component intersects $D_h$, so it has $\mu=1$. We have $\#D_h\leq 4$, so from \ref{cor:no simple curve} it follows that there are exactly two $S_0$-components in $F$, each intersecting two components of $D_h$. This eliminates the case (Y). Notice that it follows also that these two $(-1)$-curves are tips of $F$, which by \ref{lem:singular fibers} implies that $\E\cap F$ is a $(-2)$-chain between them.

Consider the case (X). We have $D_v\neq 0$, because $B^2>0$. The divisor $D_v$ is a chain and by the definition of $\sigma$ can be written as $D_v=D_0+D_1+\ldots+D_n$, where $D_0^2=-3$, $n\geq 0$ and $D_i^2=-2$ for every $1\leq i\leq n$. Let $F'$ be a fiber containing $D_v$. By \ref{cor:no simple curve} the connectedness of $D_v$ implies that each $(-1)$-curve of $F'$ intersects $D_h$. In particular, the $(-1)$-curves, and hence all components of $F'$ have $\mu=1$. It follows that $\E\cap F'=\emptyset$. We have $K\cdot D_v=1$ and $K\cdot F'=-2$, so there are exactly three $(-1)$-curves in $F'$, call them $L_2$, $L_3$ and $L_4$. We have $\sigma(F')=3$, $\sigma(F)=2$ and $\Sigma=3$ by \ref{lem:topology}(iv), so any other singular fiber has $\sigma=1$. However, the unique $(-1)$-curve of such a fiber has $\mu>1$, so cannot intersect $D_h$, hence cannot intersect $\wt D$, which is impossible. Thus $F$ and $F'$ are the only singular fibers, which implies that $\E$ is connected. Since $\mu(L_i)=1$ and $F'$ cannot contain a $(0)$-curve as a proper subdivisor, we get that one of the $L_i$'s, say $L_4$, intersects $D_n$ and two others intersect $D_0$ (it is possible that $n=0$). Each $L_i$ intersects exactly one $T_j$, so by renaming $L_i$'s we can assume that for $i=2,3,4$ we have $L_i\cdot T_i=1$. The remaining section contained in $D_h$, call it $T_1'$, is a $(-1)$-curve and intersects $D_n$. Let $M_2$ be the $(-1)$-curve of $F$ intersected by $T_4$. Denote the second $(-1)$-curve of $F$ by $M_1$. If $T_1'\cdot M_2>0$ then the contraction of $F-M_2+F'-L_4$ does not touch $T_4$ and touches $T_1'$ once, so the images of $T_4$ and $T_1'$ are disjoint sections of a $\PP^1$-ruling of a Hirzebruch surface and have self-intersections $-2$ and $0$. This is impossible, so we infer that $T_1'\cdot M_2=0$ and $T_1'\cdot M_1=1$. Now by symmetry we can assume that $T_2$ intersects $M_2$ and $T_3$ intersects $M_1$. The contraction of $F-M_1+F'-L_3$ does not touch $T_3$ and touches $T_1'$ exactly $n+1$ times. Thus as above we get a $\PP^1$-ruling of a Hirzebruch surface with two disjoint sections having self-intersections $-2$ and $n$. It follows from the properties of a Hirzebruch surface that $n=2$. Now observe that $T_4+2L_4+D_2$ and $T_3+2L_3+D_0+L_2$ are disjoint $(0)$-divisors, so they are fibers of the same $\PP^1$-ruling of $\wt S$. This contradicts the fact that $T_2$ intersects the second one and not the first one.\end{proof}

\bprop\label{prop:S' description} Let $S_0$ be the smooth locus of an exceptional singular $\Q$-homology plane $S'$. If $S'$ has at most topologically rational singularities then $\ovk(S')=\ovk(S_0)=0$ and $S'$ has a unique singular point. Moreover, either \benum[(i)]

\item $S'$ (hence $S_0$) is $\C^{**}$-ruled, its singularity is of type $A_1$ and its snc-minimal boundary $D$ is a fork with branching $(-1)$-curve and three maximal twigs: $[2]$, $[2,2,2]$ and $[2,2,2]$ (cf. \ref{ex:exception Y2c}) or

\item $S'$ (hence $S_0$) is $\C^{***}$-ruled, its singularity is of type $A_2$ and its snc-minimal boundary $D$ is a fork with branching $(-1)$-curve and three maximal twigs: $[2,2]$, $[2,2]$ and $[2,2]$. (cf. \ref{ex:exception Y1d}). \eenum\eprop

\begin{proof} We check easily that admissible chains with $d(-)$ equal to $2$, $3$ or $6$ have only one component or consist of $(-2)$-curves, so by \ref{lem:boundary description} we have only thirteen cases to consider:

\med (X0) $T_i=[2]$ for $i=1,2,3,4$ and $B^2=0$,

\med (X1) $T_i=[2]$ for $i=1,2,3,4$ and $B^2=-1$,

\med $D$ is of type (Y) with $B^2=-1$ and:

\med (Y1a) $T_1=[3], T_2=[3], T_3=[3]$,

\med (Y1b) $T_1=[3], T_2=[3], T_3=[2,2]$,

\med (Y1c) $T_1=[3], T_2=[2,2], T_3=[2,2]$,

\med (Y1d) $T_1=[2,2], T_2=[2,2], T_3=[2,2]$,

\med (Y2a) $T_1=[2], T_2=[4], T_3=[4]$,

\med (Y2b) $T_1=[2], T_2=[4], T_3=[2,2,2]$,

\med (Y2c) $T_1=[2], T_2=[2,2,2], T_3=[2,2,2]$,

\med (Y3a) $T_1=[2], T_2=[3], T_3=[6]$,

\med (Y3b) $T_1=[2], T_2=[3], T_3=[2,2,2,2,2]$,

\med (Y3c) $T_1=[2], T_2=[2,2], T_3=[6]$,

\med (Y3d) $T_1=[2], T_2=[2,2], T_3=[2,2,2,2,2]$.

\bsk Write each $T_i$ as $T_i=T_{i,1}+T_{i,2}+\ldots+T_{i,k_i}$, where $T_{i,1}$ is a tip of $D$. In cases (Y1a), (Y2a) and (Y3a) we compute $d(D)=0$, which contradicts \ref{lem:topology}(i). In each other case we specify a $\PP^1$-ruling $\pi:\ov S\to\PP^1$ with $\nu>0$ defined by some $(0)$-divisor $F_\8$ with support in $D$. By \ref{lem:topology}(iv) we have $\Sigma=\#D_h-1$. Below we list the quadruples $(F_\8,F\cdot D,\Sigma,D_v)$, where $F$ is the generic fiber and $D_v=D-\un F_\8-D_h$.

\med (X0) $(B,4,3,0)$,

\med (X1) $(T_1+2B+T_2,4,1,0)$,

\med (Y1b) $(T_1+3B+2T_{3,2}+T_{3,1},3,0,0)$,

\med (Y1c) $(T_1+3B+2T_{3,2}+T_{3,1},3,0,T_{2,1})$,

\med (Y1d) $(T_{1,2}+2B+T_{3,2},4,2,T_{2,1})$,

\med (Y2b) $(T_1+2B+T_{3,3},3,1,T_{3,1})$,

\med (Y2c) $(T_1+2B+T_{3,3},3,1,T_{3,1}+T_{2,1}+T_{2,2})$,

\med (Y3b) $(T_1+2B+T_{3,5},3,1,T_{3,1}+T_{3,2}+T_{3,3})$,

\med (Y3c) $(T_1+2B+T_{2,2},3,1,0)$,

\med (Y3d) $(T_1+2B+T_{3,5},3,1,T_{2,1}+T_{3,1}+T_{3,2}+T_{3,3})$.

\bsk Notice that $D_v$ has at most two connected components and each of them is a chain of $(-2)$-curves. Let $F$ be some singular fiber of $\pi$. The $S_0$-components of $F$ are $(-1)$-curves by \ref{lem:S_0-components are exceptional}, denote them by $L_i$, $i=1,\ldots,\sigma(F)$. We use \ref{cor:no simple curve} repeatedly.

\bcl Every $S_0$-component intersects $D_h$.\ecl

Suppose $L$ is an $S_0$-component, such that $L\cdot D_h=0$. Then $L$ intersects two $D$-components by \ref{cor:no simple curve} and these are $(-2)$-curves, so $F=[2,1,2]$. Both these $D$-components must be tips of $D$. Since $L\cdot D_h=0$ and $\nu>0$, we obtain $F\cdot D=2$, otherwise $D$ would contain a loop. This is a contradiction.

\bcl If $\mu(L_i)>1$ for some $i$ then $\sigma(F)=1$ and $\mu(L_1)=2$. \ecl

Suppose $\sigma(F)\geq 2$ and $\mu(L_1)>1$. $L_1$ intersects some $D$-component of $F$, otherwise $D_h\cdot L_1\geq 2$ and $D_h\cdot F\geq D_h\cdot (\mu(L_1)L_1+L_2)>4$, which is impossible. Thus $D_v\cap F\neq \emptyset$ and we get $4\geq D_h\cdot F\geq D_h\cdot (\mu(L_1)L_1+D_v\cap F+\mu(L_2)L_2)\geq 2+D_h\cdot (D_v\cap F)+D_h\cdot \mu(L_2)L_2$, so by (1) $\mu(L_2)=D_h\cdot L_2=1$ and $D_v\cap F$ is connected. Moreover, $\sigma(F)=2$. We get $L_2\cdot D_v>0$, because $L_2$ cannot be simple. Since $D_v$ is a $(-2)$-chain and $\mu(L_2)=1$, $L_2$ intersects $D_v$ in a tip, so $F=[1,(k),1]$ for some $k>0$ (recall that $[(k)]$ is a chain consisting of $k$ $(-2)$-curves). This contradicts $\mu(L_1)>1$.

Suppose $\sigma(F)=1$ and $\mu(L_1)>2$. Since $D_h\cdot L_1>0$, $D_h$ contains an $n$-section with $n>2$, which is possible only for (Y1b) or (Y1c). Then $FD=3$, so $D_h(\un F-L_1)=0$, hence there are no $D$-components in $F$. Thus $L_1$ is simple, a contradiction.

\bcl If $\sigma(F)>1$ then $F=[1,(k),1]$ for some $k\geq 0$. If $\sigma(F)=1$ then in cases other than (X1) $F=[2,1,2]$ and $F$ contains a $D$-component.\ecl

If $\sigma(F)>1$ then all $L_i$'s are tips of $F$ by (2). Suppose $\sigma(F)>2$. Then there are some $D$-components in $F$, otherwise $F\cdot D\geq 6$ by \ref{cor:no simple curve}. The divisor $F-\sum_iL_i$ is connected and contains a $D$-component, so there are no $\E$-components in $F$. Since $D_v$ consists of $(-2)$-curves, we get $-2=K_{\ov S}\cdot F=\sum_iK_{\ov S}\cdot L_i=-\sigma(F)$, a contradiction. Thus $\sigma(F)=2$ and both $(-1)$-curves have multiplicities one by (2), so $F=[1,(k),1]$ for some $k\geq 0$.

Assume $\sigma(F)=1$ and consider cases different from (X1). We have $\mu(L_1)=2$ by (2). There are some $D$-components in $F$, otherwise by \ref{cor:no simple curve} $L$ would meet two 2-sections contained in $D_h$, which is possible in case (X1) only. Suppose $F$ is branched. Then by \ref{rem:singular fibers} $L_1$ is a tip of $F$ and $\un F-L_1$ is one of the connected components of $D_v$, hence it must be $[2,2,2]$, which is possible for (Y3b) only. In this case $D_v$ is connected, $F\cdot D=3$ and $\Sigma=1$. In particular, there exists a fiber $F'$ with $\sigma(F')=2$ and it does not have any $D$-components, so both $S_0$-components of $F'$ meet $D_h$ at least twice, which contradicts $F\cdot D=3$. Thus $F$ is a chain, so $F=[2,1,2]$.

\bcl $\ovk(S)=0$ and $K_{\ov S}+D^\#\equiv 0$.\ecl

By (2), (3) and \ref{rem:singular fibers} every singular fiber consists of $(-1)$- and $(-2)$-curves. $\E$ is vertical, hence consist of $(-2)$-curves, so by \ref{lem:Eff-NegDef} $\ovk(S)=\ovk(S_0)=0$. The pair $(\ov S,D+\E)$ is almost minimal, so by \ref{lem:Bk properties}(v) $K_{\ov S}+D^\#+\E^\#\equiv 0$. By \ref{lem:Eff-NegDef} $K_{\ov S}+D^\#\geq_\Q 0$, so $\E=\Bk \E$ and $K_{\ov S}+D^\#\equiv 0$.

\bcl Cases other than (X0), (X1), (Y1d) and (Y2c) are impossible. $\#\E=8-B^2-\#D$.\ecl

By (4) we have $K_{\ov S}\cdot \Bk D=K_{\ov S}^2+K_{\ov S}\cdot D$, so $K_{\ov S}\cdot \Bk D\in\Z$. This excludes (Y1b), (Y1c), (Y2b), (Y3b) and (Y3c). In the remaining cases (X0), (X1), (Y1d), (Y2c) and (Y3d) the maximal twigs of $D$ are $(-2)$-chains, so by (4) $K_{\ov S}\cdot (K_{\ov S}+B)=0$. Since $\ov S$ is rational, we have $\chi(\ov S)=2+\#D+\#\E$ by \ref{lem:topology}(i) and then the Noether formula gives $12=K_{\ov S}^2+2+\#D+\#\E$, so $\#\E=8-B^2-\#D$. For (Y3d) we get $\#\E=0$, a contradiction.

\bcl $\E$ is connected. Case (X0) is impossible. \ecl

Notice that \ref{lem:KobIneq}(ii) gives $0\leq \chi(S_0)+\sum_P\frac{1}{|G_P|}\leq 1-q+\frac{q}{2}$, so if $\E$ is not connected then $q=2$ and $|G_{P_1}|=|G_{P_2}|=2$, hence $\E_1$ and $\E_2$ are $(-2)$-curves. In cases (X1) and (X0) we have $\#\E=3-B^2\geq 3$ by (5) and in case (Y2c) $\#\E=1$, so $\E$ is connected. Consider the case (Y1d). Suppose there exists a singular fiber $F$ with $\sigma(F)=1$. By (3) $F=[2,1,2]$ and there is a $D$-component in $F$, so $D_v=T_{2,1}\subseteq F$ and $F$ contains an $\E$-component. It follows that the sections $T_{1,1}$ and $T_{3,1}$ intersect $L_1$, a contradiction with $F\cdot D\leq 4$. Since $\Sigma=2$, by (3) there are only two singular fibers and they are of type $[1,(k),1]$, so $\E$ is connected because $D_v\neq 0$.

Suppose that the case (X0) occurs. Since $\Sigma=3$, there is a singular fiber $F$ with $\sigma(F)>1$, hence by (3) $F=[1,(k),1]$ for some $k\geq 0$. It is easy to see that for every such fiber $k>0$. Indeed, we know that $D_v=0$ and $L_1,L_2$ are not simple, so each is intersected by precisely two sections from $D_h$, so if $H_1,H_2\subseteq D_h$ intersect $L_1$ then $k=0$ implies that $H_1+2L_1+H_2$ gives a $\C^*$-ruling of $S_0$, a contradiction. Since $\E$ is connected, we see that there is only one fiber with $\sigma>1$. This contradicts $\Sigma=3$.

\bcl Case (X1) is impossible.\ecl

Suppose the case (X1) occurs. We have $\Sigma=1$, so by (3) there is a fiber $F_1=[1,(k),1]$ with $k\geq 0$. Suppose $k>0$. We have $D_v=0$, so $\E\subseteq F_1$ by (6) and $F_\8$ and $F_1$ are the only singular fibers. By (5) we can write $F_1=L_1+E_1+E_2+E_3+E_4+L_2$. Notice that $D_h$ consists of two $2$-sections, $T_3$ and $T_4$, and by \ref{cor:no simple curve} $D_h$ intersects $F_1-\E$ in four points. If $L_1$ intersects both $2$-sections then the contraction of $\un F_\8-T_2+F_1-L_1$ touches $T_3$ seven times, so the image of $T_3$ is a smooth 2-section on a Hirzebruch surface with self-intersection $5$, a contradiction. Thus $L_1$ intersects only one component of $D_h$, say $T_3$, hence $L_2$ intersects $T_4$. After the contraction of $\un F_\8-T_1+F_1-L_1$ the surface becomes a Hirzebruch surface and the images of the 2-sections, $T_3'$ and $T_4'$, satisfy $T_3'\cdot T_4'=2$, $T_3'^2=0$ and $T_4'^2=20$. However, $T_3'-T_4'\equiv \alpha F$ for some $\alpha\in \Z$ and a generic fiber $F$, because $T_3'$ and $T_4'$ are 2-sections. Thus $(T_3'-T_4')^2=0$, which is a contradiction. Thus $k=0$ and $\E\subseteq F_0$, where $F_0$ is a singular fiber with $\sigma(F_0)=1$. By (5) and (1) $\E$ is a $(-2)$-fork with four components. Let $M$ be the $(-1)$-curve of $F_0$. Denote the $\E$-component intersecting $M$ by $E_0$ and the branching component of $\E$ by $E_1$. Consider a new $\PP^1$-ruling of $\ov S$ given by the $(0)$-divisor $T_3+2M+T_4$. For this ruling we have $\Sigma=0$. Let $F'$ be a fiber containing $\E-E_0$. There is exactly one $(-1)$-curve $U\subseteq F'$, which is the unique $S_0$-component of $F'$. Notice that now the only possible $D$-components of $F'$ are $T_1$ and $T_2$, which are $(-2)$-curves. Since $U$ intersects some $\E$-component of $F'$, which is also a $(-2)$-curve, $U$ cannot intersect other $(-2)$-curves than $F'$, otherwise $F'=[2,1,2]$, which is not the case. We conclude that $F'$ has no $D$-components, hence $U$ intersects $E_1$ and $\mu(E_1)=\mu(U)=2$. It follows that $E_0$ intersects $F'$ only in $E_1$ and $B$ intersects $U$ in one point. Thus $U$ is a simple curve on $(\ov S,D+\E)$, a contradiction. \end{proof}

\brem\label{rem:Fujita's_Y(a,b,c)} Let us notice that smooth exceptional $\Q$-homology planes, whose description can be found in \cite[8.64]{Fujita} or \cite[4.4.4]{Miyan-OpenSurf}, can be of three types: $Y\{3,3,3\}$, $Y\{2,4,4\}$ and $Y\{2,3,6\}$, where the boundary of $Y\{d_1,d_2,d_3\}$ is a fork $$\xymatrix{{T_1}\ar@{-}[r] &{-b}\ar@{-}[r] \ar@{-}[d] &{T_2}\\ {} & {T_3} & {}} $$ with $d(T_i)=d_i$ (surfaces of type $H[-k,k]$ considered in the above references are $\C^*$-ruled). More precisely, $(b,T_1,T_2,T_3)$ is equal to $(1,[2,2],[2,2],[2,2])$ in case $Y\{3,3,3\}$, to $(0,[2],[2,2,2],[2,2,2])$ in case $Y\{2,4,4\}$ and to $(1,[2],[2,2],[2,2,2,2,2])$ in case $Y\{2,3,6\}$.

Trying to follow this notation, we see from the above proposition that singular exceptional $\Q$-homology planes are either of type $Y\{3,3,3\}$ (ruling (Y1d)) with $(b,T_1,T_2,T_3) = (1,[2,2],[2,2],[2,2])$ and $\E=[2,2]$ or of type $Y\{2,4,4\}$ (ruling (Y2c)) with $(b,T_1,T_2,T_3)= (1,[2],[2,2,2],[2,2,2])$ and $\E=[2]$. \erem

\section{Constructions}\label{s:constructions}

We now find a more precise description of rulings of type (Y2c) and (Y1d) and then use it to construct exceptional singular $\Q$-homology planes of type $Y\{2,4,4\}$ and $Y\{3,3,3\}$ respectively (cf. \ref{rem:Fujita's_Y(a,b,c)}). We produce $\Aut(\ov S,D+\E)$ -equivariant contractions $\theta\:\ov S\to \PP^2$.

\begin{figure}[h]\centering\begin{tabular}{cc}
\begin{minipage}{3in}\centering \includegraphics[scale=0.45]{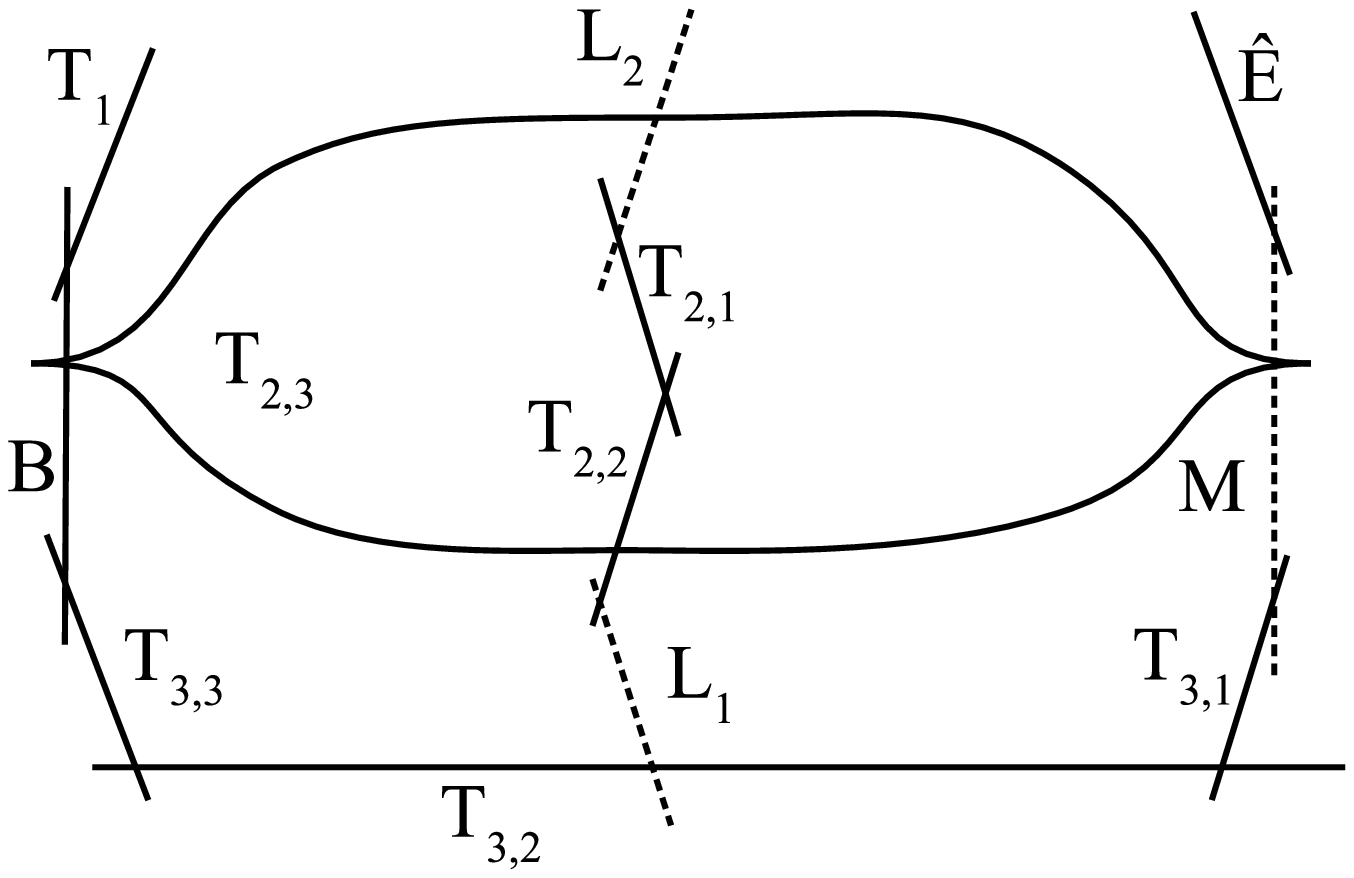}\caption{(Y2c), ruling}  \label{fig:Y2c ruling}\end{minipage}
&\begin{minipage}{3in}\centering\includegraphics[scale=0.45]{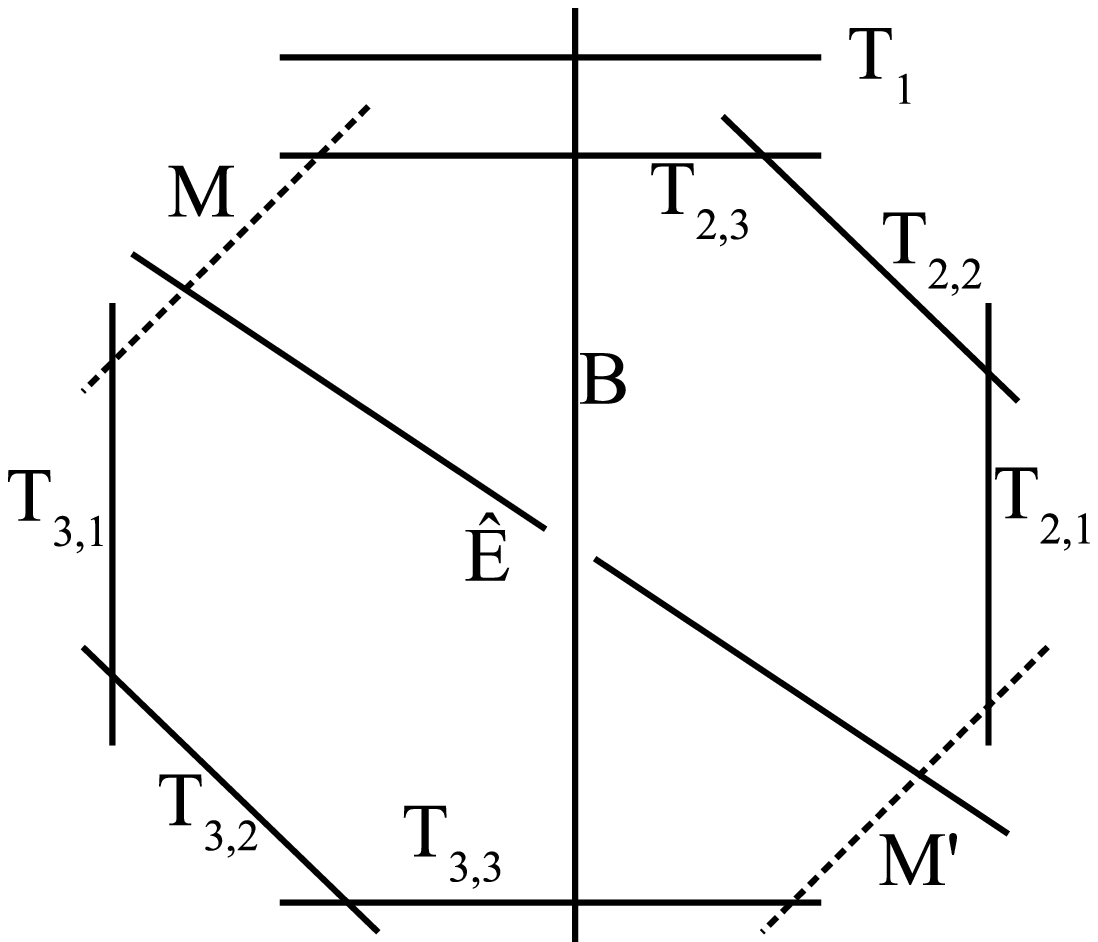}\caption{(Y2c), contraction}\label{fig:Y2c contraction}\end{minipage}
\end{tabular}\end{figure}

\begin{figure}[h]\centering\includegraphics[scale=0.45]{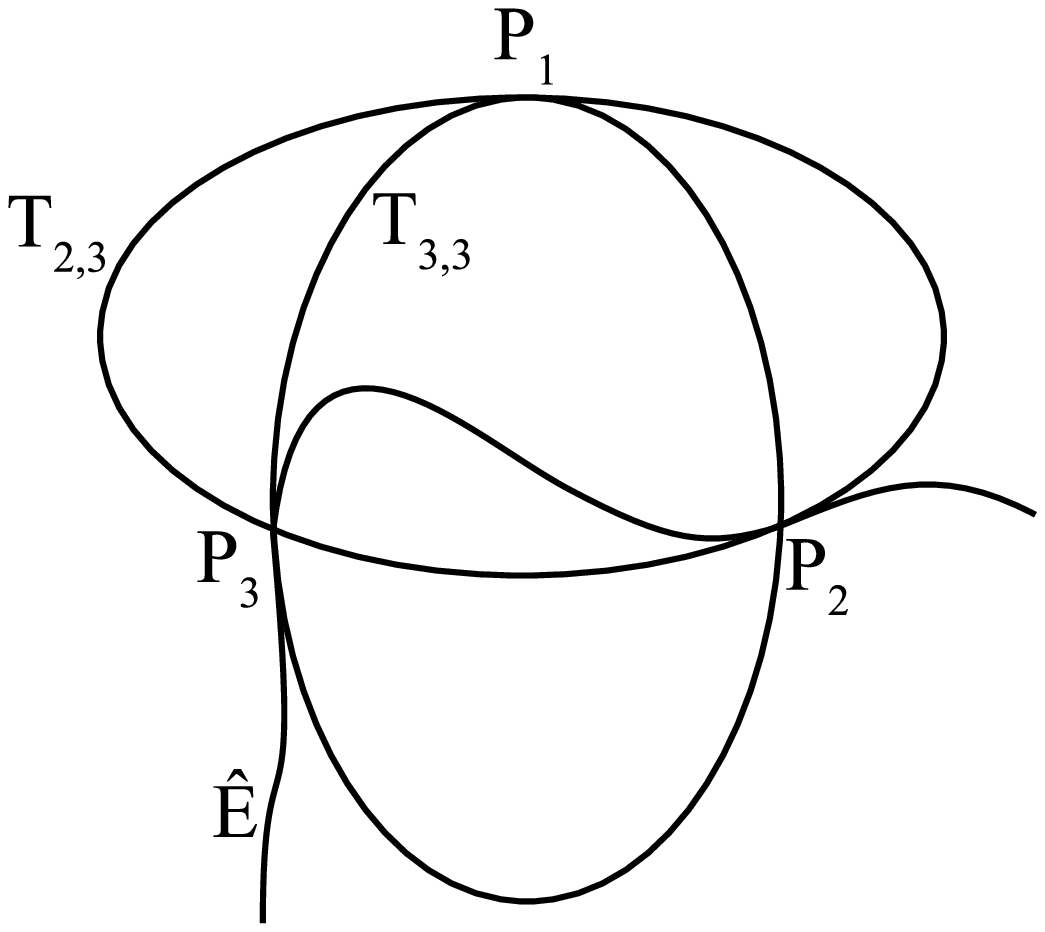}\caption{(Y2c), after contraction}\label{fig:Y2c final}\end{figure}

\blem\label{lem:Y2c description} In the case (Y2c) there are three singular fibers (see Fig.~\ref{fig:Y2c ruling}): $F_\8=T_1+2B+T_{3,3}$, $F_1=L_1+T_{2,2}+T_{2,1}+L_2$ and $F_0=T_{3,1}+M+\E$, where $\E=[2]$ and $L_1,L_2,M$ are $(-1)$-curves. The 2-section $T_{2,3}$ meets $L_2$ and $L_1\cdot T_{3,2}=1$. There is a morphism $\theta:\ov S\to \PP^1$ contracting $B+T_1+M+T_{3,1}+T_{3,2}+M'+T_{2,1}+T_{2,2}$, where $M'$ is some $(-1)$-curve, such that $\theta(T_{3,3})$ and $\theta(T_{2,3})$ are smooth conics tangent at $\theta(B)$, meeting at $\theta(T_{2,1})$ and $\theta(T_{3,1})$ and $\theta(\E)$ is a smooth conic, such that for $i=1,2$ $\theta(\E)$ intersects $\theta(T_{i,3})$ in $\theta(T_{i,1})$ with multiplicity three (see Fig.~\ref{fig:Y2c final}).\elem

\begin{proof} We use the facts showed in the proof of \ref{prop:S' description}. We have $\Sigma=1$, so by (3) there exists a fiber $F_1=[1,(k),1]$ with $k\geq 0$ and this is a unique fiber with $\sigma>1$. Since $D_v\neq 0$, $F_1$ cannot be the only singular fiber, so there exists a singular fiber $F_0$ with $\sigma(F_0)=1$. We have $F_0=[2,1,2]$ by (3). We have $\#\E=1$, $\#D_v=3$ and $D_v$ has two connected components, so $k=2$ and $F_0$ contains $\E$ and one $D$-component. Besides $F_\8$ there are no more singular fibers. Notice that $T_{2,3}$ is a 2-section intersecting the unique $(-1)$-curve of $F_0$, call it $M$, in a branching point of $\pi_{|T_{2,3}}$. Let $L_1\subset F_1$ be the $(-1)$-curve meeting $T_{2,2}$. Suppose $L_1$ meets the 2-section $T_{2,3}$ too. Then $L_2$, the second $(-1)$-curve of $F_1$, meets $T_{2,1}$ and $T_{3,2}$. The contraction of $F_\8-T_{3,3}+F_1-T_{2,2}+F_0-T_{3,1}$ touches $T_{2,3}$ five times, so the image of $T_{2,3}$ is a 2-section on a Hirzebruch surface having self-intersection $3$, which is impossible. Hence $L_1$ meets $T_{3,2}$. Let $M'$ be an exceptional component of a $\PP^1$-ruling of $\ov S$ induced by $T_1+2B+T_{2,3}$, such that $M'\cdot \E>0$. Since the structure of fibers and sections is analogous, $M'\cdot T_{2,1}=M'\cdot T_{3,3}=M\cdot \E=1$ and $M'$ does not intersect any other component of $D$ (see Fig.~\ref{fig:Y2c contraction}). Thus the chains $M+T_{3,1}+T_{3,2}$, $B+T_1$, $M'+T_{2,1}+T_{2,2}$ are disjoint and we can define $\theta:\ov S\to \PP^2$ as their contraction. Then $\theta(T_{2,3})$, $\theta(T_{3,3})$ and $\theta(\E)$ are smooth conics with prescribed properties. \end{proof}

\bcon\label{ex:exception Y2c} Let $T_{2,3}\subseteq \PP^2$ be a smooth conic and $(P_1,P_2,P_3)$ a triple of distinct points on it. This choice is unique up to an automorphism of $\PP^2$. There is a unique pair of smooth conics $(\E, T_{3,3})$, such that $P_2,P_3\in T_{3,3}\cap \E$, $T_{3,3}$ is tangent to $T_{2,3}$ at $P_1$ and $\E$ intersects $T_{i,3}$ with multiplicity three at $P_i$ for $i=2,3$ (see Fig.~\ref{fig:Y2c final}). (This can be seen as follows. Suppose $T_{2,3}=\{2yz=y^2-x^2\}$, $P_1=(0,0,1)$, $P_2=(1,-1,0)$ and $P_3=(1,1,0)$. Then the family of conics $T_{3,3}(u)$ through $P_2,P_3$ and tangent to $T_{2,3}$ at $P_1$ is one-dimensional: $T_{3,3}(u)=\{uyz=y^2-x^2\}$. The family of conics $\E(v)$ through $P_2$ and $P_3$, intersecting $T_{2,3}$ at $P_2$ with multiplicity three is one-dimensional too: $\E(v)=\{v(y^2-x^2-2yz)=z^2-yz-xz\}$. The condition for intersection at $P_3$ implies $(u,v)=(-2,\frac{1}{2})$.) We use the same names for divisors and their birational transforms. Blow up three times over $P_2$ on the intersection of $T_{2,3}$ with $\E$ and denote the subsequent exceptional curves by $T_{3,2}$, $T_{3,1}$ and $M$, similarly blow up three times over $P_3$ on the intersection of $T_{3,3}$ with $\E$ and denote the subsequent exceptional curves by $T_{2,2}$, $T_{2,1}$ and $M'$. Then blow up twice over $P_1$ so that the birational transforms of $T_{3,3}$ and $T_{2,3}$ do not meet, denote the exceptional curves by $T_1$ and $B$. Denote the resulting complete surface by $\ov S$. Define $D=T_{3,1}+T_{3,2}+T_{3,3}+T_{2,1}+T_{2,2}+T_{2,3}+T_1+B$, $S=\ov S-D$ and $S'=S/\E$. Clearly, $D$ is a fork with $\delta(D)=1$, $B^2=-1$ and other components of $D$ are $(-2)$-curves.\econ

\begin{figure}[h]\centering\begin{tabular}{cc}
\begin{minipage}{3in}\centering \includegraphics[scale=0.45]{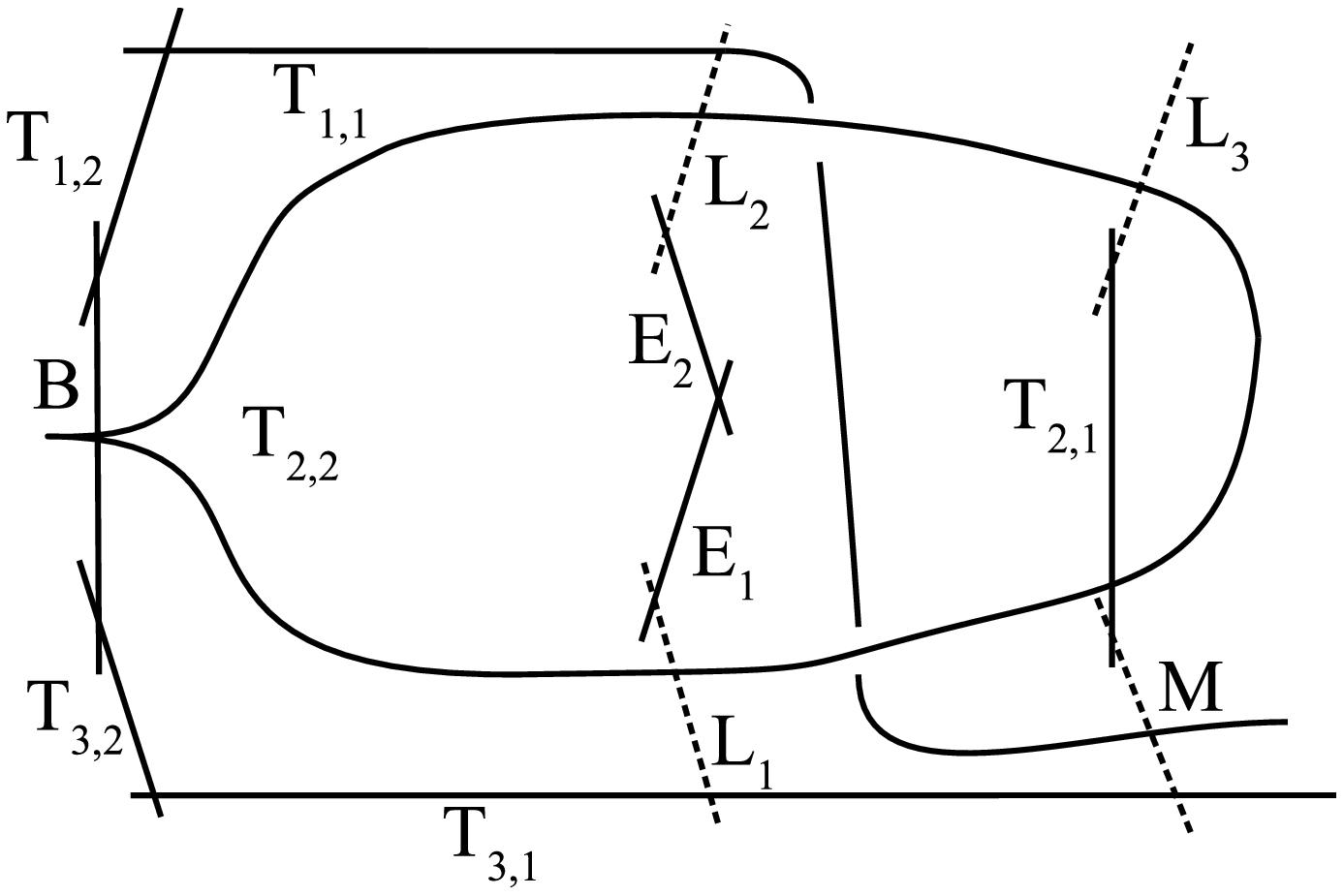}\caption{(Y1d), ruling}  \label{fig:Y1d ruling}\end{minipage}
&\begin{minipage}{3in}\centering \includegraphics[scale=0.45]{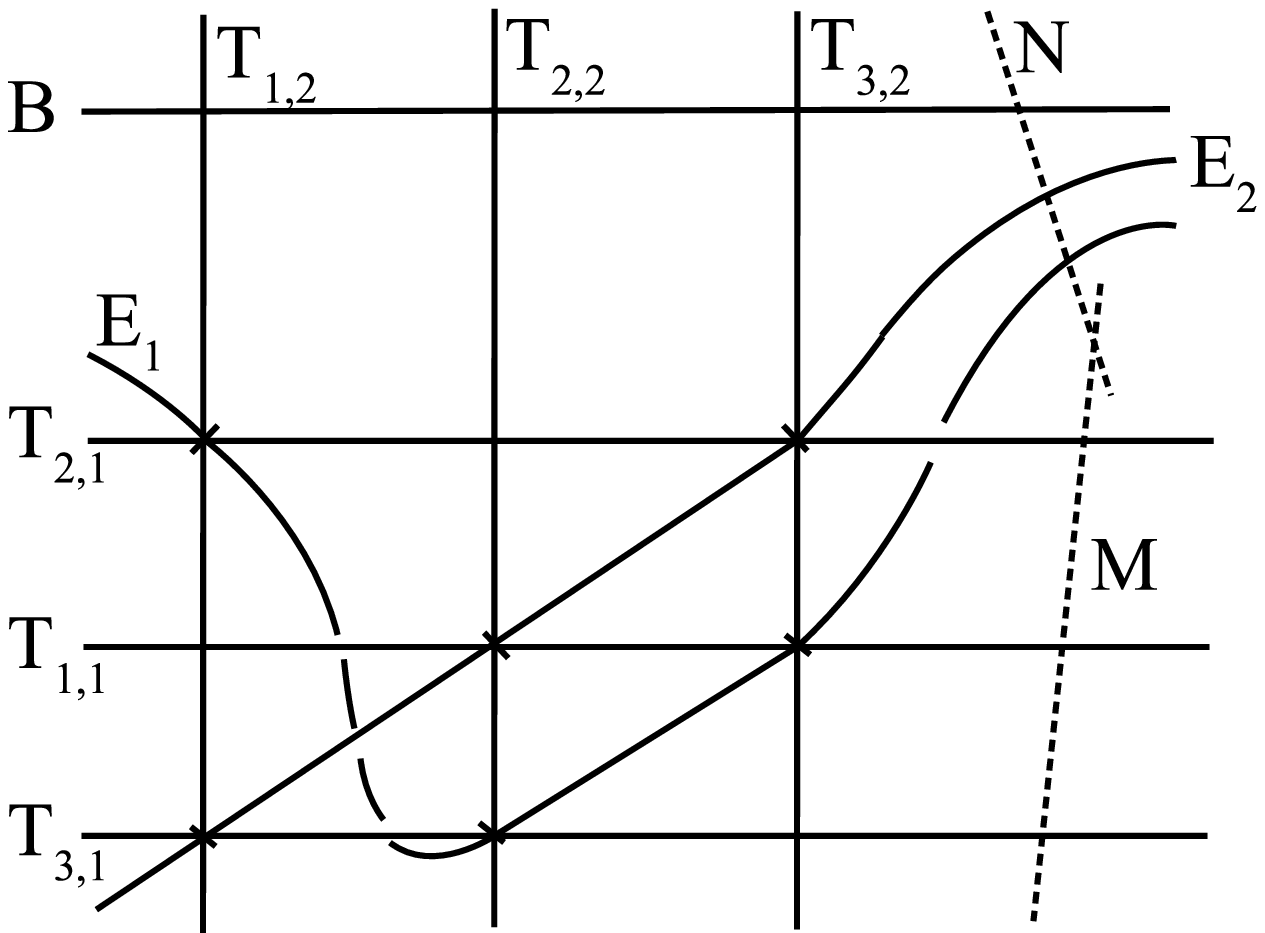}\vskip 0.2cm\caption{(Y1d), contraction}  \label{fig:Y1d other contraction}\end{minipage}
\\ \end{tabular} \end{figure}

\begin{figure}[h]\centering\includegraphics[scale=0.45]{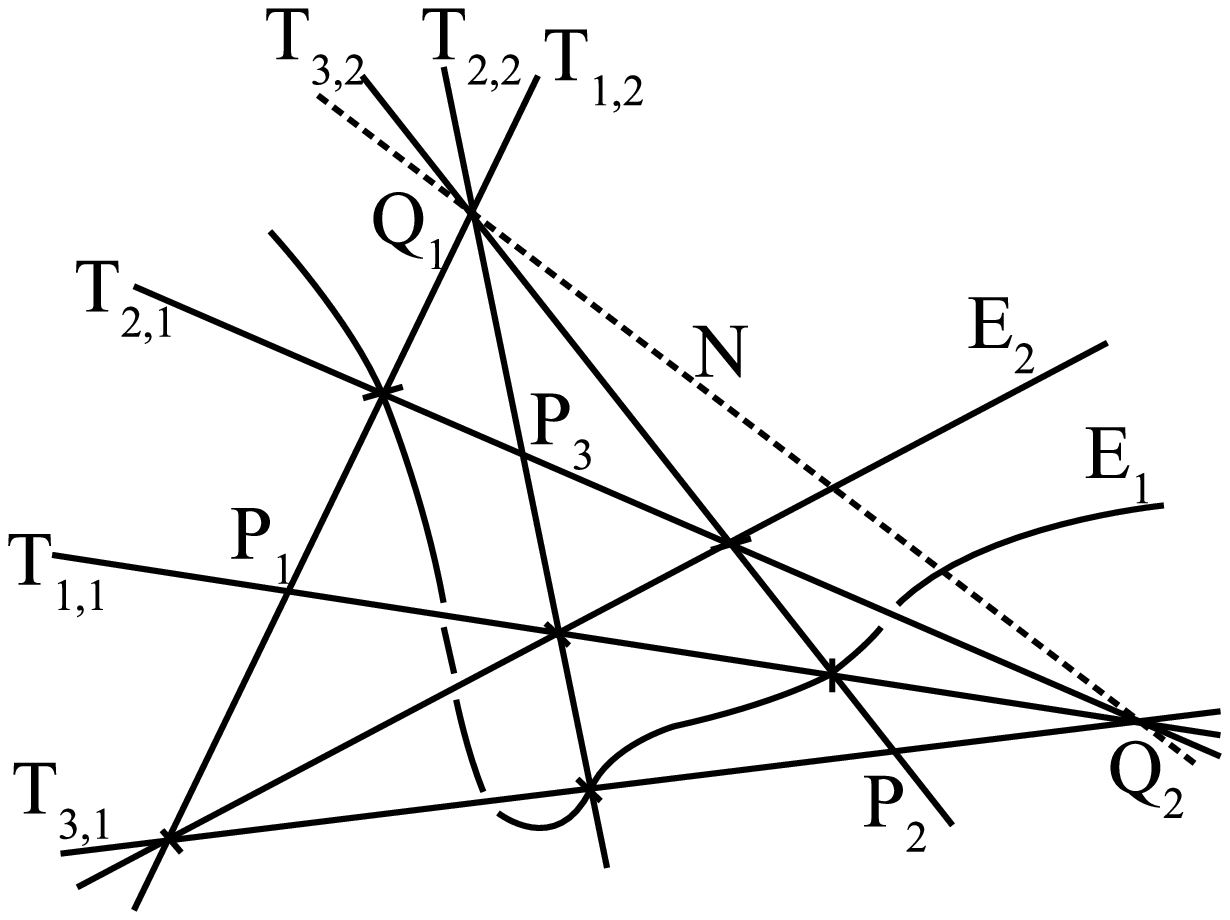}\caption{(Y1d), after contraction}  \label{fig:Y1d final}\end{figure}

\blem\label{lem:Y1d description} In the case (Y1d) there are three singular fibers (see Fig.~\ref{fig:Y1d ruling}): $F_\8=T_{1,2}+2B+T_{3,2}$, $F_1=L_1+E_1+E_2+L_2$ and $F_2=M+T_{2,1}+L_3$, where $\E=E_1+E_2=[2,2]$ and $L_1,L_2,L_3,M$ are $(-1)$-curves. $T_{3,1}\cdot M=T_{3,1}\cdot L_1=1$, $T_{1,1}\cdot L_2=T_{1,1}\cdot M=1$, $T_{2,2}\cdot L_1=T_{2,2}\cdot L_2=T_{2,2}\cdot L_3=~1$ and $T_{2,2}\cap T_{2,1}\neq T_{2,2}\cap L_3$. There exists a morphism $\theta:\ov S\to\PP^2$ contracting the divisor $B+M+L_1+L_2+L_1'+L_2'+L_1''+L_2''$ consisting of disjoint $(-1)$-curves, such that the image of $T_{1,2}+T_{2,2}+T_{3,2}$ is a triple of lines intersecting in $\theta(B)$ and the image of $T_{1,1}+T_{2,1}+T_{3,1}$ is a triple of lines intersecting in $\theta(M)$ (see Fig.~\ref{fig:Y1d final}). Moreover, $\theta(T_{1,2})\cap \theta(T_{2,1})$, $\theta(T_{2,2})\cap \theta(T_{3,1})$, $\theta(T_{3,2})\cap \theta(T_{1,1})$ lie on a line $\theta(E_1)$ and $\theta(T_{1,2})\cap \theta(T_{3,1})$, $\theta(T_{2,2})\cap \theta(T_{1,1})$, $\theta(T_{3,2})\cap \theta(T_{2,1})$ lie on a line $\theta(E_2)$.\elem

\begin{proof}

We have $\Sigma=2$, so by (3) there exist fibers $F_1=[1,(k_1),1]$, $F_2=[1,(k_2),1]$ with $k_1,k_2\geq 0$. Since $\E=[2,2]$ by (5) and singular fibers with $\sigma=1$ are of type $[2,1,2]$, we can assume that $\E\subseteq F_1$ and $k_1=2$. Since $D_v$ is connected, there are no more singular fibers besides $F_\8$, hence $T_{2,1}\subseteq F_2$ and $k_2=1$. Let $M\subseteq F_2$ be the $(-1)$-curve not intersecting $T_{2,2}$. By \ref{cor:no simple curve} $T_{1,1}+T_{3,1}$ intersects $M$, so by symmetry we can assume that $T_{3,1}$ does. Let $L_1$ be the $(-1)$-curve of $F_1$ intersecting $T_{3,1}$. The contraction of $F_\8-T_{3,2}+F_1-L_1+F_2-M$ does not touch $T_{3,1}$ and the images of $T_{3,1}$ and $T_{1,1}$ are two disjoint sections on a Hirzebruch surface, hence the image of $T_{1,1}$ must have self-intersection $2$ and we infer that the contraction touches $T_{1,1}$ exactly four times. Since $k_2=2$, it follows that $T_{1,1}$ does not intersect $L_1$ and intersects $M$ (see Fig.~\ref{fig:Y1d ruling}). Clearly, the analogous rulings of $\ov S$ induced by $F_\8'=T_{1,2}+2B+T_{2,2}$ or $F_\8''=T_{2,2}+2B+T_{3,2}$ have analogous structure of singular fibers and configuration of special horizontal components. Denote the $(-1)$-curves of the fibers of these rulings containing $\E$ as $L_1',L_2'$ and $L_1'',L_2''$ respectively. It is easy to see that $L_1,L_1',L_1'',L_2,L_2',L_2''$ are disjoint. For example, for $i=1,2$ we have $L_i\cdot F_\8'=1$, so $L_i\cdot (L_1'+L_2')=0$. Let $\omega:\ov S\to \wt S$ be the contraction of all these exceptional curves. For any $i,j,k\in\{1,2\}$ we have $\omega(T_{i,1})\cdot \omega(T_{j,2})=1$, $\omega(T_{i,j})^2=0$ and $\omega(E_k)^2=1$. We see also that $\omega(E_k)$ meets each $T_{i,j}$ once and only in points being images of curves contracted by $\omega$  (see Fig.~\ref{fig:Y1d other contraction}). Now since $b_2(\wt S)=b_2(\ov S)-6=3$, the $\PP^1$-ruling $\wt p:\wt S\to \PP^1$ induced by $\omega(T_{1,2})$ has only one singular fiber $\wt F=[1,1]$. Furthermore, $M$ is not touched by $\omega$ and $\omega(T_{1,2})\cdot M=0$, so $\wt F=M+N$, where $N$ is a birational transform of some $S_0$-component (see Fig.~\ref{fig:Y1d other contraction}). We have $\omega(T_{i,j})\cdot N=0$ and $B\cdot N=1$. If we define $\theta$ as the composition of $\omega$ with the contraction of $B+M$ then the properties of $\theta$ stated in the thesis follow (see Fig.~\ref{fig:Y1d final}).
\end{proof}

\bcon\label{ex:exception Y1d} Let $P_1=[0,1,1],P_2=[1,1,0],Q_1=[1,0,0],Q_2=[0,0,1]$ be points in $\PP^2_{(x,y,z)}$. The lines $\ov{Q_1P_1}$, $\ov{Q_1P_2}$, $\ov{Q_2P_1}$ and $\ov{Q_2P_2}$ have equations $y=z$, $z=0$, $x=0$ and $x=y$. Put $P_3=[1,\epsilon,\epsilon-1]$, where $\epsilon=-\zeta$ for some primitive third root of unity $\zeta$. Then the points $\ov{Q_1P_1}\cap \ov{Q_2P_2}=[1,1,1]$, $\ov{Q_1P_2}\cap \ov{Q_2P_3}=[\epsilon,\epsilon-1,0]$, $\ov{Q_1P_3}\cap \ov{Q_2P_1}=[0,1,\epsilon]$ lie on a line $E_2=\{(1-\epsilon)x+\epsilon y=z\}$ and the points $\ov{Q_1P_1}\cap \ov{Q_2P_3}=[1,\epsilon,\epsilon]$, $\ov{Q_1P_2}\cap \ov{Q_2P_1}=[0,1,0]$, $\ov{Q_1P_3}\cap \ov{Q_2P_2}=[1,1,\epsilon]$ lie on a line $E_1=\{z=\epsilon x\}$. Blow once in $Q_1$ and $Q_2$ and denote the exceptional curve of the first blowup by $B$. Blow once in each of the six points of intersection of lines $\ov{Q_iP_j}$ with $E_1+E_2$. Let $D$ be the divisor consisting of the proper transforms of $B$ and of lines $\ov{Q_iP_j}$. Denote the resulting surface by $\ov S$ and put $S=\ov S\setminus D$, $S'=S/\E$, where $\E=E_1+E_2$. Clearly, $D$ is a fork with $\delta(D)=1$, $B^2=-1$ and $D-B+\E$ consists of $(-2)$-curves. \econ

\brem\label{rem:configurations} Notice that the points $Q_3=E_1\cap E_2=[1,1+\epsilon,\epsilon]$, $P_1=[0,1,1]$, $P_2=[1,1,0]$ and $P_3=[1,\epsilon,\epsilon-1]$ lie on a common line $L:y=x+z$. Then the set of twelve points $\bigcup_{i=1}^3\{Q_i,P_i\}\cup (E_1\cup E_2)\cap\bigcup_{i,j=1}^3 T_{i,j}$ (with $T_{i,j}$ as on the picture \ref{fig:Y1d final}) and of nine lines $\bigcup_{i,j=1}^2\{T_{i,j}\}\cup\{E_1,E_2,L\}$ is a famous \emph{dual Hesse configuration $(12_3,9_4)$}, which is dual to the configuration of nine flexes on a smooth cubic and lines joining them (cf. \cite{Artebani-Dolgachev-Hesse_pencil} and \cite{Dolgachev-abstract_configurations}). Recall that $(a_b,c_d)$-configuration is a configuration of $a$ points and $c$ lines, such that each point lies on $b$ lines and each line contains $d$ points. This configuration has the property that each point belongs to three lines, so by the projective dual of the Sylvester-Gallai theorem, it cannot be realized in $\R\PP^2$. \erem

We now prove the theorem \ref{thm:main result}.

\begin{proof} It follows from \ref{prop:S' description} (or rather from its proof) that $S'$ is of type $Y\{2,4,4\}$ or $Y\{3,3,3\}$ (cf. \ref{rem:Fujita's_Y(a,b,c)}). If $S'$ is of type $Y\{2,4,4\}$ then the analysis of the ruling (Y2c) of $\ov S$ done in \ref{lem:Y2c description} implies that it can be constructed as in \ref{ex:exception Y2c}. The construction was determined uniquely by a choice of a smooth conic in $\PP^2$ and an ordered triple of distinct points on it, hence $S'$ of type $Y\{2,4,4\}$ is unique up to isomorphism. Clearly, the surfaces $S'$ of type $Y\{2,4,4\}$ and of type $Y\{3,3,3\}$ are non-isomorphic, because their singularities are of different type. We now prove that if $S'$ is of type $Y\{3,3,3\}$ then it can be constructed as in \ref{ex:exception Y1d}. Let $\theta:\ov S\to \PP^2$ be as in \ref{lem:Y1d description}, put $Q_1=\theta(B)$, $Q_2=\theta(M)$, $P_1=\theta(T_{1,2}\cap T_{1,1})$ and $P_2=\theta(T_{3,2}\cap T_{3,1})$, we can assume that their coordinates are as in \ref{ex:exception Y1d}. Since $P_3=\theta(T_{2,2}\cap T_{2,1})\not\in\ov{P_1Q_2}$, we can write $P_3=[1,\epsilon,u]$ for some $\epsilon,u\in\C$. The condition of collinearity of $\theta(T_{1,2})\cap \theta(T_{2,1})=[1,\epsilon,\epsilon]$, $\theta(T_{2,2})\cap \theta(T_{3,1})=[\epsilon,\epsilon,u]$, $\theta(T_{3,2})\cap \theta(T_{1,1})=[0,1,0]$ implies $u=\epsilon^2$ and the condition of collinearity of $\theta(T_{1,2})\cap \theta(T_{3,1})=[1,1,1]$, $\theta(T_{2,2})\cap \theta(T_{1,1})=[0,\epsilon,u]$, $\theta(T_{3,2})\cap \theta(T_{2,1})=[1,\epsilon,0]$ implies $\epsilon^2-\epsilon+1=0$, hence $-\epsilon$ is a primitive third root of unity. Therefore for a fixed choice of points $P_1,P_2,Q_1,Q_2$ there are two choices for $P_3$, denote them by $P_3$ and $P_3'$. The construction was determined uniquely by a choice of a quadruple of distinct points in $\PP^2$ and a primitive third root of unity, hence up to isomorphism there are at most two surfaces $S'$ of type $Y\{3,3,3\}$. For $(P_1,P_2,Q_1,Q_2)$ fixed the collinearity conditions determine the set $\{P_3,P_3'\}$. Moreover, the role of $P_1$ and $P_2$ is symmetric, so the quadruples $(P_1,P_2,Q_1,Q_2)$ and $(P_2,P_1,Q_1,Q_2)$ determine the same set $\{P_3,P_3'\}$. The automorphism $\sigma\in \Aut \PP^2$ given by \[ \begin{pmatrix} 1 & -1 & 0\\ 0 & -1 & 0\\ 0 & -1 & 1 \end{pmatrix}\] fixes $Q_1$ and $Q_2$ and changes $P_1$ with $P_2$. Since $\sigma$ changes $P_3$ with $P_3'$, we conclude that the choices of $P_3$ and $P_3'$ are equivalent.

We now check that constructions \ref{ex:exception Y1d} and \ref{ex:exception Y2c} result with singular $\Q$-homology planes with prescribed properties. In each case we have $b_1(\ov S)=0$, $b_2(\ov S)=9$ and since $d(D)\neq 0$, the components of $D+\E$ are independent in $NS(\ov S)\otimes\Q$, hence $H_2(D+\E)\to H_2(\ov S)$ is an isomorphism. The homology exact sequence of a pair $(\ov S,D)$ and the Lefschetz duality give $b_1(S)=b_3(S)=b_4(S)=0$ and $b_2(S)=\#\E$. We know that $H_2(\E)\to H_2(S)$ is a monomorphism, so the homology exact sequence of a pair $(S,\E)$ gives that $S'$ is $\Q$-acyclic. The exceptional divisors $\E$ are resolutions of singular points of type $A_1$ and $A_2$ respectively, so the constructed $S'$'s are normal. We check easily that in both cases $K_{\ov S} + D^\#$ intersects trivially with all components of $D+\E$, hence $K_{\ov S}+D^\#\equiv 0$. We check easily that in both cases $K_{\ov S} + D^\#$ intersects trivially with all components of $D+\E$, hence $K_{\ov S}+D^\#\equiv 0$ by \ref{lem:topology}(i). Then $\ovk(S)=\ovk(S_0)=0$ by \ref{lem:Eff-NegDef}.

Suppose that the smooth locus $S_0$ admits a $\C^*$-ruling. There exists a modification $(\wt S,\wt D+\wt E)\to (\ov S,D+\E)$ over $D+\E$, such that this ruling extends to a $\PP^1$-ruling $\pi:\wt S\to\PP^1$. We can assume that $\wt D+\wt E$ is $\pi$-minimal. We have $\ovk(S')\neq -\8$, so there are no sections contained in $\wt E$, hence $\wt E=\E$. The divisor $D$ does not contain components with non-negative self-intersection, which implies that this property holds for $\wt D$ too. Suppose $\#\wt D_h=1$. We have $\nu=1$ by \ref{lem:topology}(iv), so there exists a fiber $F_\8\subseteq \wt D$. Since $\wt D$ is simply connected, $F_\8$ can intersect $D_h$ only in a branching point of $\pi_{|\wt D_h}$, hence by $\pi$-minimality $F_\8=[2,1,2]$. The contractions minimalizing $\wt D$ cannot contract components of $F_\8$, hence $D$ contains two $(-2)$-tips as maximal twigs, a contradiction. Therefore we can write $\wt D_h=D_0+D_\8$ and we have $\Sigma=\nu\leq 1$ by \ref{lem:topology}(iv). If $\nu>0$ then $D_0+D_\8$ intersects the fiber contained in $\wt D$ in two different points, so this fiber is smooth by the $\pi$-minimality of $\wt D$, which contradicts the fact that all components of $\wt D$ have negative self-intersection. Thus $\Sigma=\nu=0$. Now $\ovk(S_0)=0$ implies that $F\cdot (K_{\ov S}+\wt D+\E)^-=F\cdot (K_{\wt S}+\wt D+\E)=0$, so $D_0$ and $D_\8$ are not contained in maximal twigs of $\wt D$, because are not contained in $\Supp (K_{\wt S}+\wt D+\E)^-$. The divisor $\wt D$ is simply connected, so there exists a unique fiber $F_0$, such that $F_0\cap \wt D$ is connected. By the $\pi$-minimality of $\wt D$ other singular fibers are chains intersected by $D_0$ and $D_\8$ in tips. It follows that there are at least two such fibers, otherwise $D_0$ and $D_\8$ would be contained in maximal twigs of $\wt D$. This implies that $D_0$ and $D_\8$ are branching in $\wt D$ and since exceptional components of $\wt D$ can appear only in $F_0$, after the snc-minimalization of $\wt D$ the images of $D_0$ and $D_\8$ are branching in $D$, a contradiction. \end{proof}

\bcor $\Aut Y\{3,3,3\}\cong\Z_3$ and $\Aut Y\{2,4,4\}\cong\Z_2$.\ecor

\begin{proof} Let $\eta$ be an automorphism of a surface $S'=Y\{3,3,3\}$ or $Y\{2,4,4\}$. Since $D+\E$ does not contain curves with non-negative self-intersection, $\eta_{|S_0}$ extends to $\ov \eta\in\Aut (\ov S,D+\E)$.

Suppose $S'=Y\{3,3,3\}$. We proved that $S'$ can be constructed as in \ref{ex:exception Y1d}, so we can assume that $\theta:\ov S\to \PP^2$ maps $B$ to $Q_1$ and $M$ to $Q_2$ and maps the set of nodes of $D-B$ to the fixed set of three points $\{P_1,P_2,P_3\}\subseteq \PP^2$ (we showed in the proof of the main theorem that $Q_1,Q_2,P_1,P_2$ can be fixed arbitrarily and then up to an automorphism of $\PP^2$ fixing $Q_1,Q_2$ and $\{P_1,P_2\}$ there is only one choice for $P_3$). Notice that $\ov \eta$ fixes $B$ and $M$ and acts on $\{L_1,L_1',L_2,L_2',L_3,L_3'\}$, hence descends to $\wt \eta\in\Aut \PP^2=\theta(\ov S)$ fixing $Q_1,Q_2$ and $\{P_1,P_2,P_3\}$. The automorphism of $\PP^2$ is defined uniquely by specifying the images of four points in a general position, so $\Aut S'<S_3$. However, $\sigma$ defined in the proof of \ref{thm:main result}, which fixes $Q_1,Q_2$ and exchanges $P_1$ with $P_2$, does not fix $P_3$, hence $\Aut S'<\Z_3$. We conclude that $\Aut S'\iso \Z_3$ with the generator in the coordinates as before given by $(x,y,z)\to (x-y,-\epsilon y,-\epsilon y+z)$, where $\epsilon=-\zeta$ for some primitive third root of unity $\zeta$.

Suppose $S'=Y\{2,4,4\}$. We proved that $S'$ can be constructed as in \ref{ex:exception Y2c}. Since $\ov \eta$ permutes $M$ with $M'$ and $T_{2,i}$ with $T_{3,i}$ for $i=1,2,3$, by the definition of the contraction $\theta:\ov S\to \PP^2$ it descends to $\wt \eta\in\Aut \PP^2$ fixing $P_1,\{P_2,P_3\}$ and $\{T_{2,3},T_{3,3}\}$. Notice that if $\ov \eta(T_{2,3})=T_{2,3}$ then, since $\ov \eta$ fixes $\E$ and $\E$ is tangent to $T_{2,3}$ only at $P_2$, $\ov \eta$ fixes each $P_i$, hence is an identity. It follows that if $\ov \eta$ is non-trivial then $\ov \eta(T_{2,3})=T_{3,3}$. Moreover, $\Aut S'<\Z_2$. In fact $\Aut S'\cong \Z_2$, with the generator (for conics and points as in \ref{ex:exception Y2c}) given by $(x,y,z)\to (x,-y,z)$.\end{proof}

\brem Let $M_D$ and $M$ be the 3-dimensional manifolds, which are boundaries of closures of tubular neighborhoods of $D$ and $\E$ in $S$. By \cite{Mumford} we compute that $H_1(M_D,\Z)\cong \Z_{16}\oplus\Z_2$, $H_1(M,\Z)\cong\Z_2$ for $Y\{2,4,4\}$ and $H_1(M_D,\Z)\cong \Z_9\oplus\Z_3$, $H_1(M,\Z)\cong\Z_3$ for $Y\{3,3,3\}$. Having this it in not difficult to prove that $|H_1(Y\{2,4,4\},\Z)|=4$ and $|H_1(Y\{3,3,3\},\Z)|=3$.\erem

In view of the results of \cite{DieckPetrie} it is an interesting question if the contraction $\theta:\ov S\to \PP^2$ can be chosen so that $\theta_*D+\theta_*\E$ is a sum of lines. This is clearly so for $Y\{3,3,3\}$ (cf. \ref{lem:Y2c description}) and is also possible for $Y\{2,4,4\}$. Let $\ov S$ be an snc-minimal completion of a resolution of $Y\{2,4,4\}$. We denote the twigs of $D$ as before, i.e. $T_1=[2]$, $T_2=[2,2,2]$, $T_3=[2,2,2]$. Let $\pi'\:\ov S\to \PP^1$ be a $\PP^1$-ruling induced by a $0$-curve $T_{2,3}+2B+T_{3,3}$.  Let $L_1$, $L_2$, $M$ and $M'$ be $(-1)$-curves on $\ov S$ as defined in \ref{lem:Y2c description}.

\begin{figure}[h]\centering\begin{tabular}{cc}
\begin{minipage}{3in}\centering \includegraphics[scale=0.45]{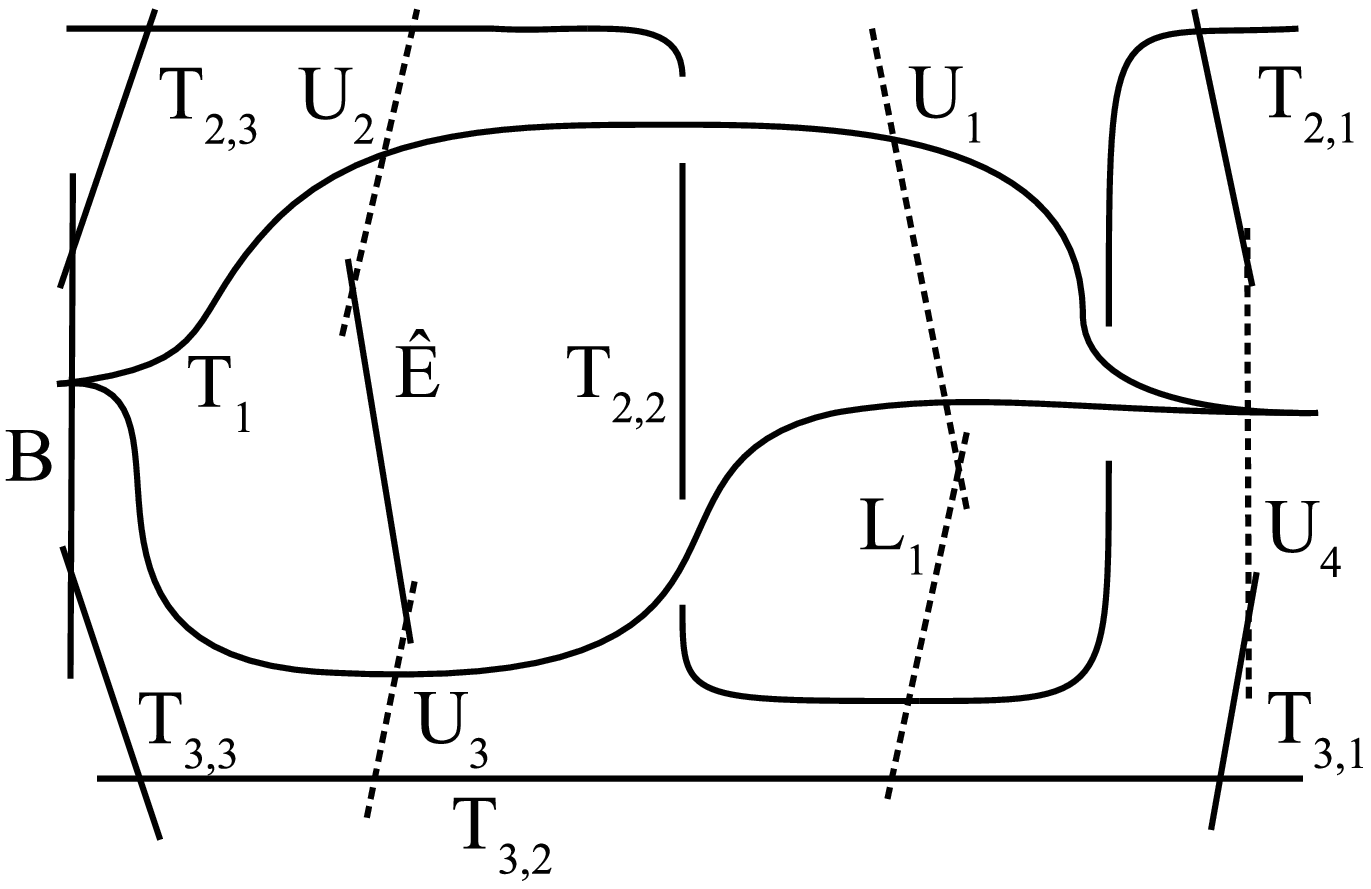}\caption{$Y\{2,4,4\}$, ruling $\pi'$}\label{fig:Y2c second ruling ruling}\end{minipage}
&\begin{minipage}{3in}\centering \includegraphics[scale=0.45]{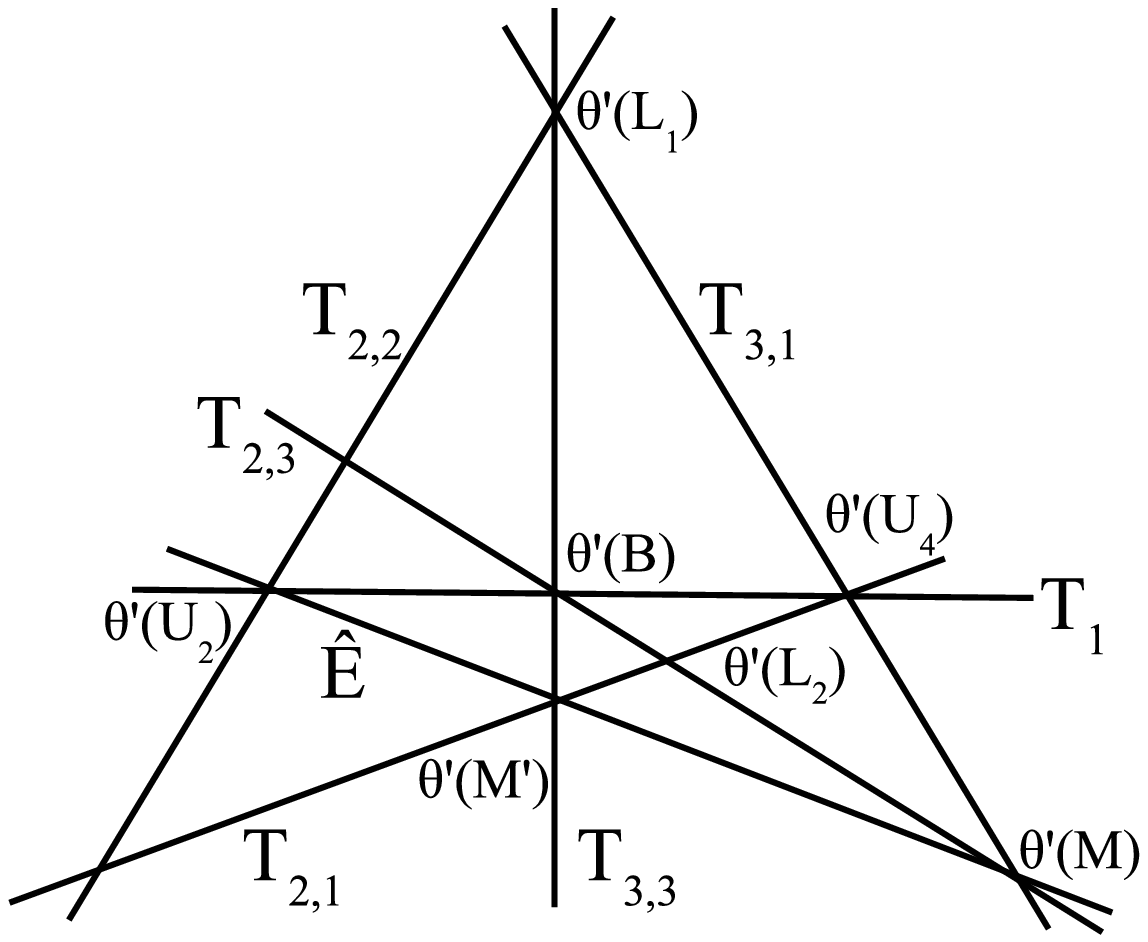}\vskip 0cm\caption{$Y\{2,4,4\}$, image of $\theta'$}  \label{fig:Y2c_second final}\end{minipage}
\\ \end{tabular} \end{figure}

\blem \label{lem:Y2c_contraction to lines}  The ruling $\pi'$ defined above has three singular fibers besides $F_\8=T_{2,3}+2B+T_{3,3}$ (see Fig.~\ref{fig:Y2c second ruling ruling}): $F_0=U_2+\E+U_3$, $F_1=U_1+L_1$ and $F_2=T_{2,1}+U_4+T_{3,1}$, where $\E=[2]$ and $U_1,U_2,U_3,U_4$ are $(-1)$-curves. We have $T_1\cdot U_2=T_1\cdot U_3=T_1\cdot U_4=1$, $T_1\cdot U_1=2$ and $T_{2,2}\cdot U_2=T_{3,2}\cdot U_3=1$. The morphism $\theta':\ov S\to \PP^2$ contracting $B+M+M'+L_1+U_1+U_4+U_2+T_{3,2}+L_2$ maps $D+\E$ into a set of lines. Namely, $\theta'(T_{2,3})$, $\theta'(T_1)$ and $\theta'(T_{3,3})$ are lines intersecting in $\theta'(B)$, $\theta'(T_{2,3})$, $\theta'(\E)$ and $\theta'(T_{3,1})$ are lines intersecting in $\theta'(M)$ and $\theta'(T_{2,1})$ is a line through $\theta'(T_{3,3})\cap \theta'(E)=\theta'(M')$ and $\theta'(T_{3,1})\cap \theta'(T_1)=\theta'(U_4)$ (see Fig.~\ref{fig:Y2c_second final}).\elem

\begin{proof} In the proof of \ref{prop:S' description} we have shown that $K_{\ov S}+D^\#\equiv 0$. Let $U$ be an $S_0$-component of some singular fiber of $\pi'$. Since $U\cdot D^\#>0$, we have $U\cdot K_{\ov S}<0$, so $U$ is a $(-1)$-curve. Then $U\cdot D^\#=1$, so computing $\Bk D$ we get $2U\cdot D_h+ U\cdot (T_{2,1}+T_{3,1})=4$. Let $F_2$ be a fiber containing $T_{2,1}$ and let $U_4$ be the $S_0$-component intersecting it. Then, since $U_4\cdot (T_{2,1}+T_{3,1})$ is even and since $F_2$ is a tree, we get $U_4\cdot T_{3,1}>0$, so in fact $U_4\cdot T_{2,1}=U_4\cdot T_{3,1}=1$. Moreover, $\un F_2=T_{2,1}+U_4+T_{3,1}$ and $U_4\cdot D_h=2$, which implies that $U_4$, having multiplicity $2$, intersects the $2$-section $T_1$. It follows that all remaining $S_0$-components $U$ have $U\cdot D_h=2$. Since $L_1\cdot \E=0$ (cf. Fig.~\ref{fig:Y2c ruling}) and the fiber $F_1$ containing $L_1$ has no $D$-components, $L_1$ intersects some $S_0$-component $U_1$, so $F_1=L_1+U_1$. We have $L_1\cdot T_{3,2}=L_1\cdot T_{2,2}=1$, so $U_1\cdot T_1=2$. The fiber $F_0$ containing $\E$ has no $D$-components, so $F_0=U_1+\E+U_2$ for some $(-1)$-curves $U_1,U_2$. By \ref{lem:topology}(iv) $\Sigma_{S_0}=2$ for $\pi'$, so there are no more singular fibers. Recall that $U_2\cdot D_h=U_3\cdot D_h=2$. It follows that each of $U_2$ and $U_3$ intersects some $1$-section contained in $D$, because if, say, $T_1\cdot U_2=0$ then the contraction of $\un F_\8-T_{3,3}+F_0-U_2+U_1+\un F_2-T_{3,1}$ does not touch $T_{3,2}$ and touches $T_{2,2}$ twice, which would result with disjoint $0$- and $(-2)$-curves as sections on a Hirzebruch surface. One can easily check that the divisors $B+L_1+U_4+U_2+T_{3,2}$ and $M+M'+L_2$ do not intersect, which implies that the contraction of $G=B+L_1+U_4+U_2+T_{3,2}+M+M'+L_2$ defines a morphism $\theta':\ov S\to \PP^2$. Each component of $D+\E$ not contained in $G$ has the intersection number with $G$ equal to three, hence each maps to a line in $\PP^2$ and the configuration of lines can be checked to be the one on the picture \ref{fig:Y2c_second final}. In particular, taking out any of the lines $\theta'(T_{2,1})$, $\theta'(T_{2,2})$ or $\theta'(T_{2,3})$ we get \emph{complete quadrangle} configurations $(4_3 6_2)$. \end{proof}

\ssk\textsl{\textsf{Acknowledgements.}} The author benefited much from discussions with his advisor dr hab. Mariusz Koras and would like to take this opportunity to express his gratitude.

\bibliographystyle{amsalpha}
\bibliography{bibl}
\end{document}